\documentclass[11pt]{article}
\usepackage{amsfonts,latexsym,amsmath,amscd,geometry}
\usepackage{amssymb}
\usepackage{latexsym}

\newcommand \nc{\newcommand}
\newtheorem{theorem}{Theorem}[section]
\newtheorem{lemma}[theorem]{Lemma}
\newtheorem{proposition}[theorem]{Proposition}
\newtheorem{corollary}[theorem]{Corollary}

\newtheorem{remark}[theorem]{Remark}

\nc{\ba}{\begin{array}}\nc{\ea}{\end{array}}
\nc{\be}{\begin{eqnarray}}\nc{\ee}{\end{eqnarray}}
\nc{\beq}{\begin{equation}}\nc{\eeq}{\end{equation}}
\nc{\bex}{\begin{eqnarray*}}\nc{\eex}{\end{eqnarray*}}
\nc{\btm}{\begin{theorem}} \nc{\etm}{\end{theorem}}
\nc{\blm}{\begin{lemma}} \nc{\elm}{\end{lemma}}
\nc{\R}{\mathbb{R}} \nc{\va}{\varepsilon} \nc{\ls}{\limits}

\def\pf{\noindent{\bf Proof.\quad}}\def\endpf{\hfill$\Box$}

\newcommand \qed {\hfill $\Box$}

\begin{document}
\title{Regularity and uniqueness of the heat flow of biharmonic maps}
\author{Jay Hineman, \ \ Tao Huang,  \ and \ Changyou Wang\\ \\
Department of Mathematics\\
University of Kentucky\\
Lexington, KY 40506, USA}
 \date{}
\maketitle

\begin{abstract}
In this paper, we first establish regularity of the heat flow of biharmonic maps into the unit sphere
$\mathbb S^L\subset\mathbb R^{L+1}$ under a smallness condition of renormalized total energy.
For the class of such solutions to the heat flow of biharmonic maps, we prove the properties
of uniqueness, convexity of hessian energy, and unique limit at $t=\infty$. We also establish both
regularity and uniqueness for the class of weak solutions $u$ to the heat flow of biharmonic maps into
any compact Riemannian manifold $N$ without boundary such that $\nabla^2 u\in L^q_tL^p_x$ for some
$p>\frac{n}2$ and $q>2$ satisfying (\ref{serrin_cond}).
\end{abstract}

\section {Introduction}
\setcounter{equation}{0}
\setcounter{theorem}{0}
For $n\geq 4$ and $L\ge k\geq 1$, let $\Omega\subset\mathbb R^n$ be a bounded smooth domain and $N\subset \mathbb R^{L+1}$ be a $k$-dimensional compact Riemannian manifold without boundary.
For $m\ge 1,\  p\ge 1$,  the Sobolev space  $W^{m,p}(\Omega,N)$ is defined by
$$W^{m,p}(\Omega, N)=\left\{v\in W^{m,p}(\Omega,\mathbb R^{L+1})\ : \ \, v(x)\in N \mbox{ for a.e. } x\in \Omega \right\}.$$
On $W^{2,2}(\Omega, N)$, there are two second order energy functionals:
$$E_2(u)=\int_{\Omega}|\Delta u|^2\ \ {\rm{and}}\ \ F_2(u)=\int_{\Omega}|(\Delta u)^{T}|^2,$$
where $(\Delta u)^T$ is the tangential component of $\Delta u$ to $T_uN$ at $u$, which is also called the tension field of $u$
(see \cite{eells-lemaire}). A map $u\in W^{2,2}(\Omega,N)$ is called an extrinsic (or intrinsic)
biharmonic map,
if $u$ is a critical point of $E_2(\cdot)$ (or $F_2(\cdot)$ respectively). It is well known that
biharmonic maps are higher-order extensions of harmonic maps, which are critical points of the
Dirichlet energy $E_1(u)=\int_\Omega |\nabla u|^2$ over $W^{1,2}(\Omega, N)$.
Recall that the Euler-Lagrange equation of (extrinsic) biharmonic maps is (see \cite{wang1} Lemma 2.1):
\begin{equation}\label{biharmonic-map}
\Delta^2 u={\mathcal N}_{\hbox{bh}}[u]:=\left[\Delta(A(u)(\nabla u,\nabla u))+2\nabla\cdot\langle\Delta u, \nabla(P(u))\rangle
-\langle\Delta(P(u)), \Delta u\rangle\right]\perp T_u N,
\end{equation}
where $P(y): \mathbb R^{L+1}\to T_y N$ is the orthogonal projection for $y\in N$, and
$A(y)(\cdot,\cdot)=\nabla P(y)(\cdot,\cdot)$ is the second fundamental form of $N$ at $y\in N$.
Throughout this paper, we use ${\mathcal N}_{\hbox{bh}}[u]$ to denote the nonlinearity in the right hand side of
the biharmonic map equation
(\ref{biharmonic-map}).

Motivated by the regularity theory of harmonic maps by Schoen-Uhlenbeck \cite{SU}, H\'elein \cite{helein},
Evans \cite{evans}, Bethuel \cite{bethuel}, Lin \cite{lin}, Rivi\`ere \cite{riviere}, and many others,  the study of biharmonic maps has attracted considerable interest and prompted a large number of interesting works by analysts during the
last several years.
The regularity of biharmonic maps to $N=\mathbb S^L$ -- the unit sphere in $\mathbb R^{L+1}$ -- was first studied by Chang-Wang-Yang \cite{chang-wang-yang}. Wang \cite{wang1, wang2, wang3} extended the main theorems of \cite{chang-wang-yang} to any compact Riemannian manifold $N$
without boundary. It asserts smoothness of biharmonic maps when the dimension $n=4$, and
the partial regularity of {\it stationary} biharmonic maps when $n\ge 5$. Here we mention in passing the
interesting works on biharmonic maps by Angelsberg \cite{angelsberg},  Strzelecki \cite{strzelecki}, Hong-Wang \cite{hong-wang},
Lamm-Rivi\`ere \cite{LR}, Struwe \cite{struwe2}, Ku \cite{ku}, Gastel-Scheven \cite{GS}, Scheven \cite{scheven1,
scheven2}, Lamm-Wang \cite{lamm-wang}, Moser \cite{moser2, moser3},
Gastel-Zorn \cite{gastel-zorn}, Hong-Yin \cite{HY}, and Gong-Lamm-Wang \cite{GLW}.

Now we describe the initial and boundary value problem for the  heat flow of biharmonic maps.
For $0<T\leq +\infty$, and $u_0\in W^{2,2}(\Omega, N)$,
a map $u\in W^{1,2}_2(\Omega\times[0,T],  N)$, i.e. $\partial_t u, \nabla^2 u \in L^2(\Omega\times [0,T])$,
is called a weak solution of the heat flow of biharmonic maps,  if $u$ satisfies in the sense of distributions
\begin{equation}\label{biharmonic-flow}
\left\{
\begin{split}
\partial_t u+\Delta^2 u =&{\mathcal N}_{\hbox{bh}}[u] \ \ \ {\rm{in}}\ \Omega\times (0, T)\\
u =& u_0 \ \ \ \ \ \ \ \ \ {\rm{on}}\ \partial_p(\Omega\times [0,T]) \\
\frac{\partial u}{\partial\nu} =& \frac{\partial u_0}{\partial\nu}\ \ \ \ \ \ \ {\rm{on}}\ \partial\Omega\times [0,T),
\end{split}
\right.
\end{equation}
where  $\nu$ denotes the outward unit normal of $\partial\Omega$. Throughout the paper,
we denote the parabolic boundary of $\Omega\times [0,T]$ by $\displaystyle\partial_p(\Omega\times [0,T])
=(\Omega\times \{0\})\cup(\partial\Omega\times (0,T))$.

The formulation of heat flow of biharmonic maps (\ref{biharmonic-flow}) remains unchanged, if
$\Omega$ is replaced by a $n$-dimensional compact Riemannian manifold $M$ with boundary $\partial M$. On the other hand,
if $\Omega$ is replaced by a $n$-dimensional compact Riemannian manifold without boundary
or a complete, non-compact Riemannian manifold without boundary $M$, then the Cauchy
problem of heat flow of biharmonic maps is considered. More precisely, if $\partial M=\emptyset$, then (\ref{biharmonic-flow}) becomes
\begin{equation}\label{biharmonic-flow2}
\left\{
\begin{split}
\partial_t u+\Delta^2 u =&{\mathcal N}_{\hbox{bh}}[u] \ {\rm{in}}\ M\times (0, T)\\
u =& u_0 \ \ \ \ \ \ \ {\rm{on}}\ M\times \{0\}.
\end{split}
\right.
\end{equation}

The Cauchy problem (\ref{biharmonic-flow2})
was first studied by Lamm \cite{lamm1}, \cite {lamm2} for $u_0\in C^\infty(M,N)$
in dimension $n=4$, where the existence of a unique, global smooth solution
is established under the condition that $\displaystyle\|u_0\|_{W^{2,2}(M)}$ is sufficiently small.
For any $u_0\in W^{2,2}(M, N)$, the existence of a unique, global weak solution of (\ref{biharmonic-flow2}), that is smooth away from finitely many times,
has been independently proved by Gastel \cite{gastel} and Wang \cite{wang4}.
We would like to point out that with suitable modifications of their proofs, the existence theorem
by \cite{gastel} and \cite{wang4} can be extended to (\ref{biharmonic-flow}) for any compact $4$-dimensional Rimannian manifold $M$ with boundary $\partial M$,
if, in additions, the trace of $u_0$ on $\partial M$ for $u_0\in W^{2,2}(M, N)$
satisfies $\displaystyle u_0|_{\partial M}\in W^{\frac72,2}(\partial M, N)$ (see \cite{huang}). Namely, there is a unique, global weak solution
$\displaystyle u\in W^{1,2}_2(M\times [0,\infty), N)$ of (\ref{biharmonic-flow}) such that
(i) $E_2(u(t))$ is monotone decreasing for $t\geq 0$; and
(ii) there exist $T_0=0<T_1<\ldots<T_k<T_{k+1}=+\infty$
such that $$\displaystyle u\in \bigcap_{i=0}^{k} C^\infty(M\times (T_i, T_{i+1}), N)\ \ {\rm{and}}\ \  \nabla u\in \bigcap_{i=0}^k
C^{\alpha}(\overline{M}\times (T_i, T_{i+1}), N),
\ \forall \ \alpha\in (0,1).$$
For dimensions $n\ge 4$,  Wang \cite{wang5} established the well-posedness of (\ref{biharmonic-flow2})
on $\mathbb R^n$ for
any $u_0:\mathbb R^n\to N$ that has sufficiently small BMO norm.
Moser \cite{moser} showed the existence of global
weak solutions $\displaystyle u\in W^{1,2}_2(\Omega\times [0,\infty), N)$  to (\ref{biharmonic-flow}) on any bounded smooth domain $\Omega\subset\mathbb R^n$ for $n\leq 8$ and $u_0\in W^{2,2}(\Omega, N)$.

Because of the critical nonlinearity in the equation (\ref{biharmonic-flow})$_1$,
the question of regularity and uniqueness for weak solutions of (\ref{biharmonic-flow}) is very challenging for dimensions $n\ge 4$.
There has been very few works in this direction.  This motivates us to
study these issues for the equation (\ref{biharmonic-flow}) in this paper.
Another motivation 
comes from our recent work \cite{huang-wang} on the heat flow of harmonic maps.
We obtain several interesting results concerning regularity, uniqueness, convexity, and unique limit at time infinity of the
equation (\ref{biharmonic-flow}), under a smallness condition of renormalized total energy.

Before stating the main theorems, we introduce some notations.

\medskip
\noindent{\bf Notations}: For $1\le p, q\le +\infty$, $0<T\le \infty$, define the Sobolev space
$$W^{1,2}_2(\Omega\times [0, T], N)
=\Big\{v\in L^2([0,T], W^{2,2}(\Omega, N)): \ \partial_t v\in L^2([0,T], L^2(\Omega))\Big\},$$
the $L^q_tL^p_x$-space
$$L^q_tL^p_x(\Omega\times [0, T], \mathbb R^{L+1})
=\Big\{f:\Omega\times [0, T]\to\mathbb R^{L+1}: \ f\in L^q([0,T], L^p(\Omega))\Big\},$$
and the Morrey space $M^{p,\lambda}_R$ for $0\le\lambda\le n+4$, $0<R\le \infty$, and $U=U_1\times U_2\subset\mathbb R^{n}\times\mathbb R$:
$${M}^{p,\lambda}_R(U)=\left\{f\in L^{p}_{\mbox{loc}}(U):
\Big\|f\Big\|_{{M}^{p,\lambda}_R(U)}<+\infty\right\},$$
where
$$\Big\|f\Big\|_{{M}^{p,\lambda}_R(U)}=\Big(\sup\limits_{(x,t)\in U}\sup\limits_{0<r<\min\{R, d(x,\partial U_1),\sqrt{t}\}}\
r^{\lambda-n-4}\int_{P_r(x,t)} |f|^p\Big)^{\frac{1}{p}},$$
and
$$B_r(x)=\{y\in\mathbb R^n: \,|y-x|\leq r\},
\ P_r(x,t)=B_r(x)\times[t-r^4,t], \ d(x, \partial U_1)=\inf_{y\in\partial U_1}|x-y|. $$
Denote $B_r$ (or $P_r$) for $B_r(0)$ (or $P_r(0)$ respectively), and
$M^{p,\lambda}(U)=M^{p,\lambda}_\infty(U)$ for $R=\infty$. We also define
the weak Morrey space $M_*^{p,\lambda}(U)$, that is the set of functions
$f$ on $U$ such that
$$\|f\|^p_{M_*^{p,\lambda}(U)}=\sup\limits_{r>0,(x,t)\in U }\,\Big\{r^{\lambda-(n+4)}
\|f\|_{L^{p,*}(P_r(x,t)\cap U)}^p \Big\}<+\infty,$$
where $L^{p,*}(P_r(x,t)\cap U)$ is the weak $L^p$-space, that is the collection of functions $v$ on $P_r(x,t)\cap U$ such that
$$
\|v\|^p_{L^{p,*}(P_r(x,t)\cap U)}=\sup\limits_{a>0}\,  \Big\{ a^p|\{z\in P_r(x,t)\cap U\ : \, |v(z)|>a\}|\Big\}<+\infty.
$$
If $N=\mathbb S^L:=\{y\in\mathbb R^{L+1}: \ |y|=1\}$, then direct calculations yield
$$\mathcal N_{\hbox{bh}}[u]=-(|\Delta u|^2+\Delta(|\nabla u|^2)+2\langle\nabla u, \nabla\Delta u\rangle)u,$$
so that for the heat flow of biharmonic maps to $\mathbb S^L$, (\ref{biharmonic-flow})$_1$ can be written into
\begin{equation}\label{biharmonic-flow1}
\partial_t u+\Delta^2 u=-(|\Delta u|^2+\Delta(|\nabla u|^2)+2\langle\nabla u, \nabla\Delta u\rangle)u.
\end{equation}

The first theorem concerns the regularity of (\ref{biharmonic-flow1}).
\begin{theorem} \label{e-regularity}
For $\frac32<p\le 2$ and $0<T<+\infty$, there exists $\epsilon_p>0$  such that if
$u\in W^{1,2}_2(\Omega\times [0,T], \mathbb S^L)$ is a weak solution of (\ref{biharmonic-flow1}) and  satisfies that,
for $z_0=(x_0, t_0)\in\Omega\times (0, T]$ and  $0<R_0\le \frac12 \min\{d(x_0,\partial\Omega), \sqrt{t_0}\}$,
\begin{equation}\label{small_norm1}
\|\nabla^2 u\|_{M^{p,2p}_{R_0}(P_{R_0}(z_0))}+\|\partial_t u\|_{M^{p,4p}_{R_0}(P_{R_0}(z_0))}
\le\epsilon_p,
\end{equation}
then $u\in C^\infty\left(P_{\frac{R_0}{16}}(z_0), \mathbb S^L\right)$, and
\begin{equation}\label{gradient_est1}
\Big|\nabla^m u(z_0)\Big|\le \frac{C\epsilon_p}{R_0^m}, \ \forall \ m\ge 1.
\end{equation}
\end{theorem}

\begin{remark}
It is an open question whether Theorem \ref{e-regularity} holds for any compact Riemannian manifold
$N$ without boundary (with $p=2$).

\end{remark}
Utilizing Theorem \ref{e-regularity},  we obtain the following uniqueness
theorem.
\begin{theorem} \label{biharmonic-flow-unique1}
For $n\ge 4$ and $\frac32<p\le 2$, there exist $\epsilon_0=\epsilon_0(p,n)>0$ and $R_0=R_0(\Omega,\epsilon_0)>0$ such that if
$u_1,u_2\in W^{1,2}_2(\Omega\times [0,T],\mathbb S^L)$ are weak
solutions of (\ref{biharmonic-flow}), with the same initial and boundary value
$u_0\in W^{2,2}(\Omega,\mathbb S^L)$,  that satisfy
\begin{equation}\label{renormal1}
\max_{i=1, 2} \Big[\|\nabla^2 u_i\|_{M^{p,2p}_{R_0}(\Omega\times (0,T))}
+\|\partial_t u_i\|_{M^{p,4p}_{R_0}(\Omega\times (0,T))}\Big]
\leq \epsilon_0,
\end{equation}
then $u_1\equiv u_2$ on $\Omega\times [0,T]$.
\end{theorem}

There are two main ingredients in the proof of Theorem \ref{biharmonic-flow-unique1}:\\
(i) The interior regularity of $u_i$ ($i=1,2$): $u_i\in C^\infty(\Omega\times (0, T), \mathbb S^L)$ and
\begin{equation}\label{grad_estimate0}
\max_{i=1, 2} |\nabla^m u_i|(x,t)\lesssim \epsilon_0
\left(\frac{1}{R_0^m}+\frac{1}{d^m(x,\partial\Omega)}+\frac{1}{t^{\frac{m}4}}\right)
\end{equation}
for any $(x,t)\in\Omega\times (0,T)$ and $m\ge 1$.\\
(ii) The energy method, with suitable applications of the Poincar\'e inequality and the second order
Hardy inequality in Lemma \ref{bi-hardy-inq-lemma} below.

\begin{remark} {\rm (i)}
We would like to point out that a novel feature of Theorem \ref{biharmonic-flow-unique1} is that the solutions may have singularities at the parabolic
boundary $\partial_p(\Omega\times [0,T])$ so that the standard argument to prove uniqueness for classical solutions
is not applicable.\\
{\rm (ii)} For $\Omega=\mathbb R^n$, if the initial data $u_0:\mathbb R^n\to N$ satisfies
that for some $R_0>0$,
$$\sup\Big\{ r^{4-n}\int_{B_r(x)}|\nabla^2 u_0|^2: \ x\in \mathbb R^n,  r\le R_0\Big\}\le\epsilon_0^2,$$
then by the local well-posedness theorem of Wang \cite{wang5} there exists $0<T_0(\approx R_0^4)$ and
a solution $u\in C^\infty(\mathbb R^n\times (0, T_0), N)$ of (\ref{biharmonic-flow2}) that
satisfies the condition (\ref{renormal1}).
\end{remark}

\smallskip
Prompted by the ideas of proof of Theorem \ref{biharmonic-flow-unique1},  we obtain the
convexity property of the $E_2$-energy along the heat flow of biharmonic maps to $\mathbb S^L$.
\begin{theorem} \label{biharmonic-flow-convex}
For $n\ge 4$, $\frac32<p\le 2$, and $1\le T\le\infty$,  there exist $\epsilon_0=\epsilon_0(p,n)>0$, $R_0=R_0(\Omega,\epsilon_0)>0$, and
$0<T_0=T_0(\epsilon_0)<T$ such that if
$u\in W^{1,2}_2(\Omega\times[0,T],\mathbb S^L)$ is a weak solution of (\ref{biharmonic-flow}),
with the initial and boundary value $u_0\in W^{2,2}(\Omega, \mathbb S^L)$,  satisfying
\begin{equation}\label{renormal2}
\|\nabla^2 u\|_{M^{p,2p}_{R_0}(\Omega\times (0,T))}+\|\partial_t u\|_{M^{p,4p}_{R_0}(\Omega\times (0,T))}
\leq \epsilon_0,
\end{equation}
then \\
(i) $E_2(u(t))$ is monotone decreasing for $t\ge T_0$;
and \\
(ii) for any $t_2\ge t_1\ge T_0$,
\begin{equation}\label{convexity2}
\int_{\Omega}|\nabla^2 ( u(t_1)-u(t_2))|^2\le C\Big[\int_{\Omega}|\Delta u(t_1)|^2-\int_\Omega |\Delta u(t_2)|^2\Big]
\end{equation}
for some $C=C(n,\epsilon_0)>0$.
\end{theorem}

A direct consequence of the convexity property of $E_2$-energy is the unique
limit at $t=\infty$ of (\ref{biharmonic-flow}).

\begin{corollary}\label{biharmonic-flow-unique2}
For $n\ge 4$ and $\frac32<p\le 2$,  there exist $\epsilon_0=\epsilon_0(p,n)>0$, and $R_0=R_0(\Omega,\epsilon_0)>0$ such that if
$u\in W^{1,2}_2(\Omega\times [0,\infty),\mathbb S^L)$ is a weak solution of (\ref{biharmonic-flow}), with the initial and
boundary value $u_0\in W^{2,2}(\Omega, \mathbb S^L)$,  satisfying the condition (\ref{renormal2}),
then there exists a biharmonic map $u_\infty\in C^\infty\cap W^{2,2}(\Omega, \mathbb S^L)$,
 with $\displaystyle(u_\infty, \frac{\partial u_\infty}{\partial\nu})=(u_0, \frac{\partial u_0}{\partial\nu})$ on $\partial \Omega$,
such that
\begin{equation}\label{unique_limit1}
\lim_{t\uparrow\infty}\|u(t)-u_\infty\|_{W^{2,2}(\Omega)}=0,
\end{equation}
and, for any compact subset $K\subset\subset \Omega$ and $m\ge 1$,
\begin{equation}\label{unique_limit2}
\lim_{t\uparrow\infty}\|u(t)-u_\infty\|_{C^m(K)}=0.
\end{equation}
\end{corollary}
\begin{remark} {\rm(i)}
We would like to remark that if  Theorem \ref{e-regularity} has been proved for any compact Riemannian manifold $N$ without boundary, then Theorem \ref{biharmonic-flow-unique1}, Theorem \ref{biharmonic-flow-convex},
and Corollary \ref{biharmonic-flow-unique2} would be true for any compact Riemannian manifold $N$
without boundary.\\
{\rm(ii)} With slight modifications of the proofs, Theorem \ref{e-regularity}, Theorem \ref{biharmonic-flow-unique1}, Theorem \ref{biharmonic-flow-convex},
and Corollary \ref{biharmonic-flow-unique2} remain to be true, if $\Omega$ is replaced by a compact Riemannian manifold
$M$ with boundary $\partial M$.\\
{\rm(iii)} If $\Omega$ is replaced by a compact or complete, non-compact Riemannian manifold $M$ with $\partial M=\emptyset$
 then Theorem \ref{e-regularity}, Theorem \ref{biharmonic-flow-unique1},
Theorem \ref{biharmonic-flow-convex},
and Corollary \ref{biharmonic-flow-unique2} remain to be true for the Cauchy problem (\ref{biharmonic-flow2}). In fact, the proof is slightly simpler than the one here,
since we don't need to use the Hardy inequalities. \\
{\rm (iv)} Schoen \cite{schoen} proved the convexity of Dirichlet energy for harmonic maps
into $N$ with nonpositive sectional curvature.  The convexity for harmonic maps into any compact manifold $N$ with small
renormalized energy was proved by \cite{huang-wang}.
In \S3 below, we will show the convexity for biharmonic maps
with small renormalized $E_2$-energy. Theorem \ref{biharmonic-flow-convex}
seems to be the first convexity result for the heat flow of biharmonic maps.\\
{\rm(v)} In general, it is a difficult question to ask whether the unique limit at $t=\infty$ holds for geometric evolution equations. Simon in his celebrated work \cite{simon} showed the unique limit at  $t=\infty$ for smooth solutions to the heat flow of harmonic maps into a real analytic manifold ($N, h$).  Corollary \ref{biharmonic-flow-unique2} seems to be first result
on the unique limit at time infinity for the heat flow of biharmonic maps.
\end{remark}

Now we consider a class of weak solutions of (\ref{biharmonic-flow}) that satisfy the smallness condition (\ref{renormal2}).
It consists of  all weak solutions $u\in W^{1,2}_2(\Omega\times [0,T], N)$ of (\ref{biharmonic-flow})
such that $\nabla^2 u\in L^q_tL^p_x(\Omega\times [0,T])$ for some
$p\ge \frac{n}2$ and $q\le \infty$ satisfying
\begin{equation}\label{serrin_cond}
\frac{n}{p}+\frac{4}{q}=2.
\end{equation}
We usually call (\ref{serrin_cond}) as Serrin's condition (see \cite{Serrin}). 
In \S 5, we will prove that if  $u$ is a weak solution of (\ref{biharmonic-flow}) such that $\nabla^2 u\in L^q_tL^p_x(\Omega\times [0,T])$
for some $p>\frac{n}2$ and $q>3$ satisfying (\ref{serrin_cond}) and $u_0\in W^{2,r}(\Omega, N)$ for some $r>\frac{n}2$,
then $u$
satisfies (\ref{renormal2}) for some $p_0>\frac32$.
Thus, for $N=\mathbb S^L$, the regularity and uniqueness for such solutions of (\ref{biharmonic-flow})  follow from Theorem
\ref{e-regularity} and Theorem \ref{biharmonic-flow-unique1}. However,
for a compact Riemannian manifold $N$ without boundary,
the regularity and uniqueness  for such a class of weak solutions of (\ref{biharmonic-flow}) require
different arguments.  More precisely, we have
\begin{theorem} \label{biharmonic-flow-serrin}
For $n\ge 4$ and $0<T\le \infty$, let $u_1,u_2\in W^{1,2}_2(\Omega\times[0,T],N)$ be weak solutions of (\ref{biharmonic-flow}),
with the same initial and boundary value $u_0\in W^{2,2}(\Omega, N)$. If, in additions, $\nabla^2 u_1,\nabla^2 u_2\in L_t^qL_x^p
(\Omega\times[0,T])$ for some $p>\frac{n}2$ and $q<\infty$ satisfying (\ref{serrin_cond}),
then $u_1, u_2\in C^\infty(\Omega\times (0, T), N)$, and $u_1\equiv u_2$ in $\Omega\times[0,T]$.
\end{theorem}
\begin{remark} {\rm(i)} It is a very interesting question to ask whether Theorem \ref{biharmonic-flow-serrin}
holds  for the end-point case
$p=\frac{n}2$ and $q=\infty$.\\
{\rm(ii)} If $u_0\in W^{2,r}(\Omega,N)$ for some $r>\frac{n}2$, then the local existence of solutions $u$
of (\ref{biharmonic-flow}) such that $\nabla^2 u\in L^q_tL^p_x(\Omega\times [0,T])$ for some
$p>\frac{n}2$ and $q<\infty$ satisfying (\ref{serrin_cond}) can be shown by the fixed point argument
similar to \cite{FJR} \S 4. We leave it to interested readers. 
\end{remark}

For dimension $n=4$, by applying Theorem \ref{e-regularity-serrin} (with $p=2$ ($=\frac{n}2$) and $q=\infty$) and
the second half of the proof of Theorem \ref{biharmonic-flow-unique1}, we obtain the following uniqueness result.
\begin{corollary}\label{4d_unique}
For $n=4$ and $0<T\le\infty$, there exists $\epsilon_1>0$ such that if $u_1
\ {\rm{and}}\ u_2\in W^{1,2}_2(\Omega\times [0,T], N)$ are weak solutions of (\ref{biharmonic-flow}),
under the same initial and boundary value $u_0\in W^{2,2}(\Omega, N)$,  satisfying
\begin{equation}\label{no_jump}
\limsup_{t\downarrow t_0^+} E_2(u_i(t))\le E_2(u_i(t_0))+\epsilon_1, \ \forall\  t_0\in [0, T),
\end{equation}
for $i=1,2$. Then $u_1\equiv u_2$ in $\Omega\times [0,T)$. In particular, the uniqueness holds among weak solutions of (\ref{biharmonic-flow}),
whose $E_2$-energy is monotone decreasing for $t\ge 0$.
\end{corollary}

We would like to point out that for the Cauchy problem (\ref{biharmonic-flow2})
of heat flow of biharmonic maps on a compact
$4$-dimensional Riemannian manifold $M$ without boundary, Corollary \ref{4d_unique} has been recently proven by Rupflin \cite{rupflin} through a different
argument.

Concerning the convexity and unique limit of (\ref{biharmonic-flow}) at $t=\infty$ in dimension $n=4$, we have
\begin{corollary}\label{4d_convex} For $n=4$, there exist $\epsilon_2>0$ and $T_1>0$
such that if $u\in W^{1,2}_2(\Omega\times (0,+\infty), N)$ is a weak solution
of (\ref{biharmonic-flow}), with the initial-boundary value $u_0\in W^{2,2}(\Omega, N)$, satisfying
\begin{equation}
E_2(u(t))\le \epsilon_2^2, \ \forall\ t\ge 0, \label{mono_decrease}
\end{equation}
then
(i) $E_2(u(t))$ is monotone decreasing for $t\geq T_1$;\\
(ii) for $t_2\ge t_1\ge T_2$, it holds
$$\int_\Omega |\nabla^2(u(t_1)-u(t_2))|^2\le C\left(E_2(u(t_1))-E_2(u(t_2))\right)$$
for some $C=C(\epsilon_2)>0$; and\\
(iii) there exists a biharmonic map $u_\infty\in C^\infty\cap W^{2,2}(\Omega, N)$, with
$\displaystyle (u_\infty, \frac{\partial u_\infty}{\partial\nu})=(u_0, \frac{\partial u_0}{\partial\nu})$ on $\partial\Omega$, such
that $\displaystyle \lim_{t\rightarrow\infty}\|u(t)-u_\infty\|_{W^{2,2}(\Omega)}=0$, and for any $m\ge 1$, $K\subset\subset\Omega$,
$\displaystyle \lim_{t\rightarrow\infty}\|u(t)-u_\infty\|_{C^m(K)}=0$.
\end{corollary}

It is easy to see that the condition (\ref{mono_decrease}) holds for any solution $u\in W^{1,2}_2(\Omega\times [0,\infty), N)$ of (\ref{biharmonic-flow})  that satisfies $E_2(u(t))\le E_2(u_0)$ for $t\ge 0$ (e.g., the solution by \cite{gastel}
and \cite{wang4}) and $E_2(u_0)\le\epsilon_2^2$.

The paper is written as follows. In \S2, we will prove the $\epsilon$-regularity Theorem \ref{e-regularity}
for weak solutions of (\ref{biharmonic-flow}) under the assumption (\ref{small_norm1}). In \S3, we will show
both convexity and uniqueness property for biharmonic maps with small $E_2$-energy. In \S4, we  will
prove the uniqueness Theorem \ref{biharmonic-flow-unique1}, the convexity Theorem \ref{biharmonic-flow-convex},
and the unique limit Theorem \ref{biharmonic-flow-unique2}.  In \S5, we will discuss weak solutions $u$
of (\ref{biharmonic-flow}) such that $\nabla^2 u\in L^q_tL^p_x(\Omega\times [0,T])$ for some
$p\ge\frac{n}2$ and $q\ge 2$ satisfying (\ref{serrin_cond}), and prove Theorem \ref{biharmonic-flow-serrin}, Corollary \ref{4d_unique},  and Corollary \ref{4d_convex}. In \S6 Appendix, we will sketch a
proof for higher-order regularities of the heat flow of biharmonic maps.

\section {$\epsilon$-regularity}
\setcounter{equation}{0}
\setcounter{theorem}{0}
This section is devoted to the proof of Theorem \ref{e-regularity}, i.e., the regularity of heat flow of biharmonic maps to $\mathbb S^L$
under the smallness condition (\ref{small_norm1}). The idea is motivated by \cite{chang-wang-yang} on the regularity
of stationary biharmonic maps to $\mathbb S^L$.

The first step is to rewrite (\ref{biharmonic-flow1}) into the form where nonlinear terms
are of divergence structures, analogous to \cite{chang-wang-yang} on the equation of biharmonic maps to
$\mathbb S^L$.
As in \cite{chang-wang-yang}, we divide the nonlinearities in
(\ref{biharmonic-flow1}) into four different types: for $1\le\alpha\le L+1$,
\begin{equation}\label{e-reg-eq1}
\begin{split}
&T_{11}^\alpha=\left(u_j^\alpha\Delta u^\beta(u^\beta-c^\beta)\right)_j
\ {\rm{or}}\ \left(u_j^\beta\Delta u^\alpha(u^\beta-c^\beta)\right)_j, \ T_{12}^\alpha
=\left((u^\alpha-c^\alpha)u_i^\beta u^\beta_{ij}\right)_j,\\
&T_{21}^\alpha=\Delta\left((u^\alpha-c^\alpha)|\nabla u|^2\right),
\ T_{22}=\Delta\left((u^\beta-c^\beta)\Delta u^\beta\right),\\
&T_{23}^\alpha=\Delta\left(u^\alpha(u^\beta-c^\beta)\Delta u^\beta\right)
\ {\rm{or}}\  \Delta\left(u^\beta (u^\beta-c^\beta)\Delta u^\alpha\right),\\
& T_{33}=\left((u^\beta-c^\beta)u^\beta_j\right)_{jii},\ T_{41}^\alpha=\left(u^\alpha \partial_t u^\beta-u^\beta \partial_t u^\alpha\right)\left(u^\beta-c^\beta\right),
\end{split}
\end{equation}
where the upper index $\alpha,\ \beta$ denotes the component of a vector, the lower index $i,\ j$
denotes the differentiation in the direction $x_i,\ x_j$, $c^\alpha\in \mathbb R^{L+1}$ is a constant, and the
Einstein convention of summation is used.

\begin{lemma}\label{bi-div-structure} The equation 
(\ref{biharmonic-flow1}) is equivalent to
\begin{eqnarray}
\label{e-reg-eq2}
\partial_t u^\alpha+\Delta^2 u^\alpha=\mathcal F_\alpha(T_{11}^\alpha,\
T_{12}^\alpha, \ T_{21}^\alpha, \ T_{22}, \ T_{23}^\alpha, \ T_{33}, \ T_{41}^\alpha),
\ 1\le\alpha\le L+1, 
\end{eqnarray}
where $\mathcal F_\alpha$ denotes a linear function of its arguments such that
the coefficients can be bounded independent of $u$.

\end{lemma}

\pf We follow \cite{chang-wang-yang} Proposition 1.2 closely.
First, by Lemma 1.3 of \cite{chang-wang-yang}, we have that,  for every fixed $\alpha$,
\begin{equation}\label{e-reg-eq3}
c^\alpha\Delta\left(|\nabla u|^2\right)\mbox{ and } \left(u^\alpha_j|\nabla u|^2\right)_j\mbox{ are linear functions of }\
 T_{11}^\alpha, T_{12}^\alpha, T_{21}^\alpha,T_{22}, T_{23}^\alpha,  T_{33},
\end{equation}
whose coefficients can be bounded independent of $u$.
For $1\le\alpha\le L+1$, set
\begin{equation}\label{e-reg-eq4}
S_1^\alpha=u^\alpha|\Delta u|^2,\ S_2^\alpha=2u^\alpha u^\beta_j\left(\Delta u^\beta\right)_j,\ S_3^\alpha=u^\alpha\Delta\left(|\nabla u|^2\right).
\end{equation}
Differentiation of $|u|=1$ gives
\begin{equation}\label{e-reg-eq5}
u^\gamma u^\gamma_j=0,\ u^\gamma\Delta u^\gamma+|\nabla u|^2=0.
\end{equation}
By the equation (\ref{biharmonic-flow}), we have
\begin{equation}\label{e-reg-eq6}
u^\alpha\Delta^2 u^\beta+u^\alpha \partial_t u^\beta=u^\beta\Delta^2u^\alpha+u^\beta \partial_t u^\alpha,
\ 1\le\alpha,\beta\le L+1.
\end{equation}
It follows from (\ref{e-reg-eq5}) and (\ref{e-reg-eq6}) that
\begin{equation}\label{e-reg-eq7}
\begin{split}
\frac{S_2^\alpha}{2}=&u^\alpha u^\beta_j (\Delta u^\beta)_j
=u^\beta_j\left(u^\alpha (\Delta u^\beta)_j-u^\beta\left(\Delta u^\alpha\right)_j\right)\\
=&u^\beta_j\left(u^\alpha (\Delta u^\beta)_j-u^\beta\left(\Delta u^\alpha\right)_j-u^\alpha_j\Delta u^\beta+u^\beta_j\Delta u^\alpha\right)+u^\beta_j\left(u^\alpha_j\Delta u^\beta-u^\beta_j\Delta u^\alpha\right)\\
=&\left\{\left(u^\beta-c^\beta\right)\left(u^\alpha (\Delta u^\beta)_j-u^\beta\left(\Delta u^\alpha\right)_j-u^\alpha_j\Delta u^\beta+u^\beta_j\Delta u^\alpha\right)\right\}_j\\
&+\left(u^\beta-c^\beta\right)\left(u^\alpha \partial_t u^\beta-u^\beta \partial_t u^\alpha\right)+u^\beta_j\left(u^\alpha_j\Delta u^\beta-u^\beta_j\Delta u^\alpha\right)\\
=&\left\{\left(u^\beta-c^\beta\right)\left(u^\alpha\Delta u^\beta-u^\beta\Delta u^\alpha\right)\right\}_{jj}
-\left\{u^\beta_j\left(u^\alpha\Delta u^\beta-u^\beta\Delta u^\alpha\right)\right\}_j\\
&-2\left\{\left(u^\beta-c^\beta\right)\left(u^\alpha_j\Delta u^\beta-u^\beta_j\Delta u^\alpha\right)\right\}_j
+u_j^\beta\left(u^\alpha_j\Delta u^\beta-u^\beta_j\Delta u^\alpha\right)+T_{41}^\alpha\\
=&-\left\{u^\beta_j\left(u^\alpha\Delta u^\beta-u^\beta\Delta u^\alpha\right)\right\}_j+u_j^\beta\left(u^\alpha_j\Delta u^\beta-u^\beta_j\Delta u^\alpha\right)+G_\alpha(T_{11}^\alpha, T_{21}^\alpha, T_{23}^\alpha, T_{41}^\alpha),
\end{split}
\end{equation}
where $G_\alpha$ is a linear function of its arguments whose coefficients can be bounded independent of $u$.
By (\ref{e-reg-eq3}) and (\ref{e-reg-eq5}), we have
\begin{equation}\label{e-reg-eq8}
\begin{split}
S_3^\alpha=&\left(u^\alpha-c^\alpha\right)\Delta\left(|\nabla u|^2\right)+c^\alpha\Delta\left(|\nabla u|^2\right)\\
=&\Delta\left(\left(u^\alpha-c^\alpha\right)|\nabla u|^2\right)-2\left(u^\alpha_j(|\nabla u|^2\right)_j
-\Delta u^\alpha u^\beta \Delta u^\beta+c^\alpha\Delta\left(|\nabla u|^2\right)\\
=&
-\Delta u^\alpha u^\beta\Delta u^\beta+H_\alpha(T_{11}^\alpha, T_{12}^\alpha, T_{21}^\alpha,
T_{22}, T_{23}^\alpha, T_{33}),
\end{split}
\end{equation}
where $H_\alpha$ is a linear function of its arguments whose coefficients can be bounded independent of $u$.
By (\ref{e-reg-eq8}), the definition of $S_1^\alpha$, and (\ref{e-reg-eq7}), we have
\begin{equation}\label{e-reg-eq9}
\begin{split}
S_1^\alpha+S_3^\alpha=&\left(u^\alpha\Delta u^\beta-u^\beta\Delta u^\alpha\right)\Delta u^\beta+
\mathcal H_\alpha(T_{11}^\alpha, T_{12}^\alpha, T_{21}^\alpha,
T_{22}, T_{23}^\alpha, T_{31}),\\
=&\left\{\left(u^\alpha\Delta u^\beta-u^\beta\Delta u^\alpha\right)u_j^\beta\right\}_j
-\left(u^\alpha_j\Delta u^\beta-u^\beta_j\Delta u^\alpha\right)u_j^\beta\\
&-\left(u^\alpha\Delta u_j^\beta-u^\beta\Delta u_j^\alpha\right)u_j^\beta+\mathcal H_\alpha(T_{11}^\alpha, T_{12}^\alpha, T_{21}^\alpha,
T_{22}, T_{23}^\alpha, T_{31}),\\
=&-\frac{S_2^\alpha}{2}-\frac{S_2^\alpha}{2}+L_\alpha(T_{11}^\alpha, T_{12}^\alpha, T_{21}^\alpha,
T_{22}, T_{23}^\alpha, T_{33}, T_{41}^\alpha),
\end{split}
\end{equation}
where $L_\alpha$ is a linear function of its arguments whose coefficients can be bounded independent of $u$.
Therefore we obtain
$$S_1^\alpha+S_2^\alpha+S_3^\alpha=L_\alpha(T_{11}^\alpha, T_{12}^\alpha, T_{21}^\alpha,
T_{22}, T_{23}^\alpha, T_{33}, T_{41}^\alpha).$$
This completes the proof.
\endpf

\medskip
Next we recall some basic properties of the heat kernel for $\Delta^2$ in $\mathbb R^n$,  and the definition of
Riesz potentials on $\mathbb R^{n+1}$, and the definition of BMO space and John-Nirenberg's inequality (see \cite{john-nirenberg}).
Let $b(x,t)$  be the fundamental solution of
$$(\partial_t +\Delta^2) v=0 \ {\rm{in}}\ \mathbb R^{n+1}_+.$$
Then we have (see \cite{koch-lamm} \S 2.2):
$$b(x,t)=t^{-\frac{n}{4}}g\left(\frac{x}{t^{\frac{1}{4}}}\right),
\  {\rm{with}}\ g(\xi)=(2\pi)^{-\frac{n}{2}}\int_{\mathbb R^n} e^{i\xi \eta-|\eta|^4}, \ \xi\in\mathbb R^n, $$
and the estimate
\begin{equation}\label{e-reg-eq15}
\Big|\nabla^m b(x,t)\Big|\leq C\left(|t|^{\frac{1}{4}}+|x|\right)^{-n-m}, \ \forall\ (x,t)\in\mathbb R^{n+1}_+, \, \forall \ m\geq 1.
\end{equation}
We equip $\mathbb R^{n+1}$ with the parabolic distance $\delta$:
$${\bf\delta}((x,t),(y,s))=|t-s|^{\frac{1}{4}}+|x-y|, \ (x,t), \ (y,s)\in\mathbb R^{n+1}.$$
For $0\le \alpha\le n+4$, define the Riesz potential of order $\alpha$ on ($\mathbb R^{n+1}$, $\delta$) by
\begin{equation}\label{e-reg-eq18}
I_{\alpha}(f)(x,t)=\int_{\mathbb R^{n+1}}\left(|t-s|^{\frac{1}{4}}+|x-y|\right)^{\alpha-n-4}|f|(y,s),
\ (x,t)\in\mathbb R^{n+1}.
\end{equation}
For any open set $U\subset\mathbb R^{n+1}$, let BMO($U$) denote the space of functions of bounded mean oscillations:
$f\in {\rm{BMO}}(U)$ if
\begin{equation}
\label{e-reg-eq12}
[f]_{\mbox{BMO}(U)}:=\sup\Big\{-\!\!\!\!\!\!\int_{P_r(z)}|f-f_{P_r(z)}|: \ P_r(z)\subset U\Big\}<+\infty,
\end{equation}
where $\displaystyle-\!\!\!\!\!\!\int_{P_r(z)}=\frac{1}{|P_r(z)|}\int_{P_r(z)}$ and $\displaystyle f_{P_{r}(z)}=-\!\!\!\!\!\!\int_{P_{r}(z)}f$
denotes the average of $f$ over $P_r(z)$. By the celebrated John-Nirenberg inequality (see \cite{john-nirenberg}),
we have that if $f\in {\rm{BMO}}(U)$, then for any $1<q<+\infty$  it holds
\begin{equation}\label{JN-ineq}
\sup\Big\{\left(-\!\!\!\!\!\!\int_{P_{r}(z)}|f-f_{P_{r}(z)}|^q\right)^{\frac{1}{q}}:\ \ P_r(z)\subset U\Big\}\leq
C(q)\Big[f\Big]_{\rm{BMO}(U)}.
\end{equation}

Now we are ready to prove the $\epsilon$-regularity  for the heat flow of biharmonic maps
to $\mathbb S^L$.

\begin{proposition} \label{e-regular-biharmonic-flow}{\it For any $\frac32<p\le 2$, there exists $\epsilon_p>0$
such that if $u:P_4\rightarrow\mathbb S^L$ is a weak solution of
(\ref{biharmonic-flow1})
and satisfies
\begin{equation}
\label{e-condition}
\sup\limits_{(x,t)\in P_3, 0<r\leq 1} r^{2p-n-4}\int_{P_r(x,t)}\left(|\nabla^2 u|^{p}+r^{2p}|\partial_t u|^p\right)\leq\epsilon_p^p,
\end{equation}
then $u\in C^{\infty}(P_{\frac{1}{2}},\mathbb S^L)$, and
\begin{equation}\label{gradient_estimate}
\Big\|\nabla^m u\Big\|_{C^0(P_\frac12)}\le C(p, n, m), \ \forall\ m\ge 1.
\end{equation}
}
\end{proposition}

\pf We first establish  H\"older continuity of $u$ in $P_\frac34$.
It is based on the decay estimate.\\
\noindent{\it Claim}. There exist $\epsilon_p>0$ and $\theta_0\in(0,\frac{1}{2})$ such that
\begin{equation}\label{e-reg-decay}
\Big[u\Big]_{\rm{BMO}(P_{\theta_0})}\leq
\frac{1}{2} \Big[u\Big]_{\rm{BMO}(P_{2})}.
\end{equation}
In order to establish (\ref{e-reg-decay}), we first want to prove that there exists $q>1$ such that
\begin{equation} \label{bmo-decay}
-\!\!\!\!\!\!\int_{P_{\theta r}(z_0)}|u-u_{P_{\theta r}(z_0)}|\le C\left(\theta^{-(n+4)}\epsilon_p+\theta\right)
\left(-\!\!\!\!\!\!\int_{P_{r}(z_0)}|u-u_{P_{r}(z_0)}|^q\right)^{\frac1{q}}
\end{equation}
holds for any $0<\theta\le\frac12$,  $z_0\in P_1$, and $0<r\le 2$.

 By translation and scaling, it suffices to show (\ref{bmo-decay}) for $z_0=(0,0)$ and $r=2$.
First, we need to extend $u$ from $P_1$ to $\mathbb R^{n+1}$. Let the extension, still denoted by $u$, be
such that
$$|u|\le 2 \ {\rm{in}}\ \mathbb R^{n+1}, \ \ u=0 \ {\rm{outisde}}\ P_2,$$
and
$$\int_{\mathbb R^{n+1}} |\nabla^2 u|^p+|\partial_t u|^p
\lesssim \int_{P_2} |\nabla^2 u|^p+|\partial_t u|^p.$$
For $1\le \alpha\le L+1$, let  $w_{ij}^\alpha:\mathbb R^{n+1}_+\to\mathbb R$ 
be solutions of
\begin{equation}\label{e-reg-eq10}
\partial_t w_{ij}^\alpha+\Delta^2 w_{ij}^\alpha=T_{ij}^\alpha \ \ \mbox{ in }\ \mathbb R^{n+1}_+;\
w_{ij}^\alpha=0\ \   \mbox{ on }\ \mathbb R^n\times\{0\}
\end{equation}
for $ij\in\{11, 12, 21, 23, 41\}$, and and $w_{kk}:\mathbb R^{n+1}_+\to\mathbb R$ be solutions of
\begin{equation}\label{e-reg-eq10}
\partial_t w_{kk}+\Delta^2 w_{kk}=T_{kk}\ \ \mbox{ in }\ \mathbb R^{n+1}_+;\
w_{kk}=0\ \ \mbox{ on }\ \mathbb R^n\times\{0\}
\end{equation}
for $k\in\{2, 3\}$.
Define $v:P_1\to\mathbb R^{L+1}$ by letting
 $$v^\alpha=u^\alpha-\mathcal F_\alpha(w_{11}^\alpha, w_{12}^\alpha, w_{21}^\alpha, w_{22}, w_{23}^\alpha, w_{33}, w_{41}^\alpha),
\ 1\le\alpha\le L+1.$$
Here $\mathcal F_\alpha$ is the linear function given by Lemma \ref{bi-div-structure}.
By (\ref{e-reg-eq2}), we have
\begin{equation}\label{e-reg-eq11}
\partial_t v+\Delta^2 v=0\ \mbox{ in }\ P_1.
\end{equation}
It follows from (\ref{e-reg-eq10}) and the Duhamel formula that  for $1\le \alpha\le L+1$,
\begin{equation}\label{e-reg-eq16}
\begin{cases}
w_{ij}^\alpha (x,t)=\int_{\mathbb R^{n}\times [0,t]}b(x-y,t-s) T_{ij}^\alpha(y,s), \ & ij\in\{11, 12, 21, 23, 41\},\\
w_{kk}(x,t)=\int_{\mathbb R^{n}\times [0,t]}b(x-y,t-s) T_{kk}(y,s), \ & k\in\{2, 3\}.
\end{cases}
\end{equation}
Set $c^\alpha=u_{P_2}^\alpha$ in (\ref{e-reg-eq1}). Then it is easy to see $|c^\alpha|\le 1$.
Now we can estimate $w_{12}^\alpha$ by ($w_{11}^\alpha$ can be estimated similarly):
\begin{equation}\label{e-reg-eq17}
\begin{split}
|w_{12}^\alpha(x,t)|=&\left|\int_{\mathbb R^{n}\times [0,t]}\nabla_jb(x-y,t-s)(u^\alpha-u_{P_2}^\alpha)u_i^\beta u^\beta_{ij}(y,s)\right|\\
\lesssim &\int_{\mathbb R^{n+1}}\left(|t-s|^{\frac{1}{4}}+|x-y|\right)^{-n-1}|u-u_{P_2}||\nabla u||\nabla^2u|(y,s)\\
\lesssim &I_3\left(\chi_{P_2}|u-u_{P_2}||\nabla u||\nabla^2u|\right)(x,t),
\end{split}
\end{equation}
where $\chi_{P_2}$ is the characteristic function of $P_2$.

By the estimate of Riesz potentials in $L^{q}$-spaces (see also \S5 below),
we have that for any $f\in L^{q}$, $1<q<+\infty$,  $I_\alpha(f)\in L^{\tilde{q}}$,
where $\frac{1}{\tilde{q}}=\frac{1}{q}-\frac{\alpha}{n+4}.$
As $p>\frac32$, we can check that for sufficiently large $q_1>1$,
there exists $\widetilde{q_1}>1$ such that
$$\frac{1}{\widetilde {q_1}}=\frac{1}{p}+\frac{1}{2p}+\frac{1}{q_1}-\frac{3}{n+4}.$$
Hence we obtain
\begin{equation}\label{e-reg-eq19}
\begin{split}
\Big\|w_{12}^\alpha\Big\|_{L^{\widetilde{q_1}}(P_2)}\leq C\Big\|u-u_{P_2}\Big\|_{L^{q_1}(P_2)}
\Big\|\nabla u\Big\|_{L^{2p}(P_2)}\Big\|\nabla^2 u\Big\|_{L^{p}(P_2)}
\leq C\epsilon_p\Big\|u-u_{P_2}\Big\|_{L^{q_1}(P_2)}.
\end{split}
\end{equation}
Next we can estimate $w_{21}^\alpha$ by
($w_{22}$ and $w_{23}^\alpha$ can be estimated similarly):
\begin{equation}\label{e-reg-eq20}
\begin{split}
|w_{21}^\alpha(x,t)|=&\left|\int_{\mathbb R^n\times [0,t]}\Delta b(x-y,t-s)(u^\alpha-u_{P_2}^\alpha)|\nabla u|^2(y,s)\right|\\
\lesssim &\int_{\mathbb R^{n+1}}\left(|t-s|^{\frac{1}{4}}+|x-y|\right)^{-n-2}|u-u_{P_2}||\nabla u|^2(y,s)\\
\lesssim &I_2\left(\chi_{P_2}|u-u_{P_2}||\nabla u|^2\right)(x,t).
\end{split}
\end{equation}
For $q_2>1$ sufficiently large, there exists $\widetilde{q}_2>1$ be such that
$$\frac{1}{\widetilde{q}_2}=\frac{1}{p}+\frac{1}{q_2}-\frac{2}{n+4}.$$
Hence we obtain
\begin{equation}\label{e-reg-eq21}
\begin{split}
\Big\|w_{21}^\alpha\Big\|_{L^{\widetilde{q}_2}(P_2)}\leq
C\Big\|u-u_{P_2}\Big\|_{L^{q_2}(P_2)}
\Big\||\nabla u|^2\Big\|_{L^p(P_2)}
\leq C\epsilon_p\Big\|u-u_{P_2}\Big\|_{L^{q_2}(P_2)}.
\end{split}
\end{equation}
For $w_{33}$, we have
\begin{equation}\label{e-reg-eq22}
\begin{split}
|w_{33}(x,t)|=&\left|\int_{\mathbb R^n\times [0,t]}\Delta b_j(x-y,t-s)(u^\beta-u_{P_2}^\beta)u^\beta_j(y,s)\right|\\
\lesssim &\int_{\mathbb R^{n+1}}\left(|t-s|^{\frac{1}{4}}+|x-y|\right)^{-n-3}|u-u_{P_2}||\nabla u|(y,s)\\
\lesssim &I_1\left(\chi_{P_2}|u-u_{P_2}||\nabla u|\right).
\end{split}
\end{equation}
For $q_3>1$ sufficiently large, there exists $\widetilde{q}_3>1$ such that
$$\frac{1}{\widetilde{q}_3}=\frac{1}{2p}+\frac{1}{q_3}-\frac{1}{n+4}.$$
Hence we obtain
\begin{equation}\label{e-reg-eq23}
\begin{split}
\Big\|w_{33}\Big\|_{L^{\widetilde{q}_3}(P_2)}
\leq C\Big\|u-u_{P_2}\Big\|_{L^{q_3}(P_2)}
\Big\|\nabla u\Big\|_{L^{2p}(P_2)}
\leq C\epsilon_p\Big\|u-u_{P_2}\Big\|_{L^{q_3}(P_2)}.
\end{split}
\end{equation}
For $w_{41}^\alpha$, we have
\begin{equation}\label{e-reg-eq24}
\partial_t w_{41}^\alpha+\Delta^2 w_{41}^\alpha=\left(u^\alpha \partial_t u^\beta-u^\beta\partial_t u^\alpha\right)\left(u^\beta-u_{P_2}^\beta\right).
\end{equation}
By the Duhamel formular, we have
$$w_{41}^\alpha(x,t)
=\sum_\beta\int_0^t\int_{\mathbb R^n} b(x-y,t-s) \left(u^\alpha\partial_t u^\beta-u^\beta \partial_t u^\alpha\right)\left(u^\beta-u_{P_2}^\beta\right)(y,s),$$
so that by applying the Young inequality we obtain
\begin{eqnarray}\label{w-estimate}
\|w_{41}\|_{L^{\widetilde{q_4}}(\mathbb R^n\times [0,2])}
&\lesssim&\|b\|_{L^1(\mathbb R^n\times [0,2])}\left(\sum_{\alpha, \beta} \Big\| (u^\alpha \partial_t u^\beta-u^\beta \partial_t u^\alpha)(u^\beta-u_{P_2}^\beta)\Big\|_{L^{\widetilde{q_4}}(\mathbb R^n\times [0,2])}\right)\nonumber\\
&\lesssim& \|\partial_t u\|_{L^p(P_2)}\|u-u_{P_2}\|_{L^{q_4}(P_2)},
\end{eqnarray}
where $q_4>\frac{p}{p-1}$ and $1<\widetilde{q_4}<p$ satisfy
$$\frac{1}{\widetilde{q_4}}=\frac{1}{p}+\frac{1}{q_4}.$$
Set
$$q=\max\left\{q_1,q_2,q_3, q_4\right\}>1\mbox{ and }
\widetilde q=\min\left\{\widetilde{q}_1,\widetilde{q}_2,\widetilde{q}_3, \widetilde{q_4}\right\}>1.$$
By (\ref{e-reg-eq19}), (\ref{e-reg-eq21}), (\ref{e-reg-eq23}) and (\ref{w-estimate}), we have
\begin{equation}\label{e-reg-eq28}
\sum_{ij=11,12,21,23,41}^{1\le\alpha\le L+1}\|w_{ij}^\alpha\|_{L^{\widetilde q}(P_2)}
+\sum_{k=2}^3\|w_{kk}\|_{L^{\widetilde q}(P_2)}
\leq C\epsilon_p\Big\|u-u_{P_2}\Big\|_{L^{q}(P_2)}.
\end{equation}
On the other hand, by the standard estimate on $v$, we have that for any
$0<\theta<1$,
\begin{equation}\label{e-reg-eq29}
\left(-\!\!\!\!\!\!\int_{P_{\theta}}|v-v_{P_{\theta}}|^{\widetilde q}\right)^{\frac{1}{\widetilde q}}
\leq
C\theta\left(-\!\!\!\!\!\!\int_{P_{1}}|v-v_{P_{1}}|^{q}\right)^{\frac{1}{q}}
\leq
C\theta\Big\|u-u_{P_2}\Big\|_{L^{q}(P_2)}.
\end{equation}
Adding (\ref{e-reg-eq28}) and (\ref{e-reg-eq29}) together and applying the H\"older inequality, we obtain
\begin{equation}\label{e-reg-eq30}
-\!\!\!\!\!\!\int_{P_{\theta}}|u-u_{P_{\theta}}|\leq
\left(-\!\!\!\!\!\!\int_{P_{\theta}}|u-u_{P_{\theta}}|^{\widetilde q}\right)^{\frac{1}{\widetilde q}}\leq
C\left(\theta^{-(n+4)}\epsilon_p+\theta\right)
\left(-\!\!\!\!\!\!\int_{P_{2}}|u-u_{P_{2}}|^q\right)^{\frac{1}{q}}.
\end{equation}
This implies (\ref{bmo-decay}).

Now we indicate how (\ref{e-reg-decay}) follows from (\ref{bmo-decay}). It follows from the Poincar\'e inequality and
(\ref{e-condition}) that $u\in {\rm{BMO}}(P_3)$, and hence by (\ref{JN-ineq}) we have
\begin{equation} \label{bmo-decay1}
-\!\!\!\!\!\!\int_{P_{\theta r}(z_0)}|u-u_{P_{\theta r}(z_0)}|\le C\left(\theta^{-(n+4)}\epsilon_p+\theta\right)
\Big[u\Big]_{\rm{BMO}(P_2)}
\end{equation}
holds for any $0<\theta\le\frac12$,  $z_0\in P_1$, and $0<r\le 1$.
Taking supremum of (\ref{bmo-decay1}) over all $z_0\in P_\theta$ and $0<r\le 1$, we obtain
\begin{equation}\label{bmo-decay2}
\Big[u\Big]_{\rm{BMO}(P_\theta)}\le
C\left(\theta^{-(n+4)}\epsilon_p+\theta\right)
\Big[u\Big]_{\rm{BMO}(P_2)}.
\end{equation}
If we choose $\theta=\theta_0\in (0,\frac12)$ and $\epsilon_p$ small enough so that
$$C\left(\theta_0^{-(n+4)}\epsilon_p+\theta_0\right)\leq \frac{1}{2},$$
then (\ref{bmo-decay2}) implies (\ref{e-reg-decay}).

It is standard that iterating (\ref{e-reg-decay}) yields the H\"older continuity of $u$
by using the Campanato theory \cite{campanato}. The higher-order regularity then follows from
the hole-filling type argument and the bootstrap argument, which will be sketched
in Proposition 6.1 of \S6 Appendix. After this, we have that
$u\in C^\infty(P_\frac12, \mathbb S^L)$ and the estimate (\ref{gradient_estimate}) holds.
\endpf

\medskip
\noindent{\bf Proof of Theorem \ref{e-regularity}}. By the definition of Morrey spaces, for $z_0=(x_0,t_0)\in \Omega\times (0, T)$ and
$R_0\le \frac12\min\{d(x_0,\partial\Omega), \sqrt{t_0}\}$,
we have
\begin{equation}\label{renormal20}
\sup_{z\in P_{\frac{R_0}2}(z_0), \ r\le \frac{R_0}2} r^{2p-(n+4)}\int_{P_r(z)}(|\nabla^2 u|^p+r^{2p}|\partial_t u|^p)
\le\epsilon_p^p.
\end{equation}
Consider $v(x,t)=u(x_0+\frac{R_0}8 x, t_0+(\frac{R_0}8)^4 t): P_4\to\mathbb S^L$. It is easy to check that
$v$ is a weak solution of (\ref{biharmonic-flow1}) and satisfies  (\ref{e-condition}). Hence Proposition 2.2 implies
that $v\in C^\infty(P_\frac12, \mathbb S^L)$ and satisfies (\ref{gradient_estimate}). After rescaling, we see
that $u\in C^\infty(P_{\frac{R_0}{16}}(z_0),\mathbb S^L)$ and the estimate (\ref{gradient_est1}) holds. \qed

\medskip
Since biharmonic maps are steady solutions of the heat flow of biharmonic maps,
as a direct consequence of Theorem \ref{e-regularity} we have the following
$\epsilon$-regularity for biharmonic maps to $\mathbb S^L$.

\begin{corollary} \label{e-regular-biharmonic-map}{\it For $\frac32<p\le 2$, there exist
$\epsilon_p>0$ and $r_0>0$ such that if $u\in W^{2,p}(\Omega, \mathbb S^L)$
is a weak solution of (\ref{biharmonic-map})
and satisfies
\begin{equation}
\label{e-condition-bi}
\sup\limits_{x\in \Omega}\sup\limits_{0<r\leq \min\{r_0, d(x,\partial\Omega)\}}
r^{2p-n}\int_{B_r(x)}|\nabla^2 u|^p\leq\epsilon_p^p,
\end{equation}
then $u\in C^{\infty}(\Omega,\mathbb S^L)$, and
\begin{equation}\label{gradient_estimate1}
|\nabla^m u(x)|\le C\epsilon_p\Big(\frac{1}{r_0^m}+\frac{1}{d^m(x,\partial\Omega)}\Big),
\ \forall \ m\ge 1.
\end{equation}
}
\end{corollary}

\begin{remark}  For $p=2$, Corollary \ref{e-regular-biharmonic-map}  was first proved by
Chang-Wang-Yang \cite{chang-wang-yang}. For biharmonic maps into any
compact Riemannian manifold $N$ without boundary,
Corollary \ref{e-regular-biharmonic-map} was proved by \cite{wang1, wang3} for $p=2$.

\end{remark}

\section {Convexity and uniqueness of biharmonic maps}
\setcounter{equation}{0}
\setcounter{theorem}{0}

We will show the convexity and
uniqueness properties for biharmonic maps with small energy, which are
the second-order extensions of the theorems on harmonic maps with small
energy by Struwe \cite{struwe}, Moser \cite{moser1}, and Huang-Wang \cite{huang-wang}.

Consider the Dirichlet  problem for a biharmonic map $u\in W^{2,2}(\Omega, N)$:
\begin{equation}\label{biharmonic-map-bd}
\left\{
\begin{split}
\Delta^2 u=&\mathcal N_{\hbox{bh}}[u]
\quad\quad\ \mbox{in }\ \ \Omega\\
\Big(u, \frac{\partial u}{\partial\nu}\Big)
=&\Big(u_0, \frac{\partial u_0}{\partial\nu}\Big) \quad\mbox{on } \ \ \partial \Omega.
\end{split}
\right.
\end{equation}
where $u_0\in W^{2, 2}(\Omega, N)$ is given.

We recall the second order
Hardy inequality.
\begin{lemma}\label{bi-hardy-inq-lemma}{\it  There is $C>0$ depending
only on $n$ and $\Omega$ such that  if $f\in {W}^{2,2}_0(\Omega)$, then
\begin{equation}\label{hardy-ineq}
\int_{\Omega}\frac{|f(x)|^2}{d^4(x,\partial\Omega)}
\leq C\int_{\Omega}|\nabla^2 f(x)|^2.
\end{equation}
}
\end{lemma}

\pf For simplicity, we indicate a proof for the case $\Omega=B_1$ -- the unit ball
in $\mathbb R^n$. The readers can refer to \cite{ER} for a proof of general domains.
By approximation, we may assume
$f\in C^{\infty}_0(B_1)$. Writing the left hand side of (\ref{hardy-ineq}) in
 spherical coordinates, integrating over $B_1$, and using the H$\ddot{\mbox{o}}$lder inequality, we obtain
\begin{equation}\label{unq-bi-map-eq1}
\begin{split}
\int_{B_1}\frac{|f(x)|^2}{(1-|x|)^4}=&\int_0^1\int_{\mathbb S^{n-1}}\frac{|f|^2(r,\theta)}{(1-r)^4}r^{n-1}\,dH^{n-1}(\theta) dr\\
=&-\int_0^1\int_{\mathbb S^{n-1}}\frac{1}{3(1-r)^3}\left(2ff_rr^{n-1}+|f|^2(n-1)r^{n-2}\right)\,dH^{n-1}(\theta)dr\\
\leq& -\int_0^1\int_{\mathbb S^{n-1}}\frac{2}{3(1-r)^3}ff_rr^{n-1}\,dH^{n-1}(\theta)dr\\
\leq& C\int_0^1\int_{\mathbb S^{n-1}}\frac{|f||f_r|r^{n-1}}{(1-r)^3}\,dH^{n-1}(\theta)dr\\
\leq &C\int_{B_1}\frac{|f(x)||\nabla f(x)|}{(1-|x|)^3}\\
\leq& C\left(\int_{B_1}\frac{|f(x)|^2}{(1-|x|)^4}\right)^{\frac{1}{2}}
\left(\int_{B_1}\frac{|\nabla f(x)|^2}{(1-|x|)^2}\right)^{\frac{1}{2}}.
\end{split}
\end{equation}
Thus, by using the first-order Hardy inequality, we obtain
\begin{equation}\label{unq-bi-map-eq2}
\begin{split}
\int_{B_1}\frac{|f(x)|^2}{(1-|x|)^4}\leq C\int_{B_1}\frac{|\nabla f(x)|^2}{(1-|x|)^2}
\leq C\int_{B_1}|\nabla^2 f(x)|^2.
\end{split}
\end{equation}
This yields (\ref{hardy-ineq}). \endpf

\medskip
Now we introduce the Morrey spaces in $\mathbb R^n$. For $1\le l<+\infty$, $0<\lambda\le n$,
and $0<R\le +\infty$, $f\in M^{l,\lambda}_R(\Omega)$ if $f\in L^l_{\hbox{loc}}(\Omega)$
satisfies
$$
\|f\|_{M^{l,\lambda}_R(\Omega)}^l
:=\sup_{x\in \Omega}\sup_{0<r\le \min\{R, d(x,\partial \Omega)\}}
\Big\{r^{\lambda-n}\int_{B_r(x)}|f|^l\Big\}<+\infty.$$

We have the convexity property of biharmonic maps with small energy.
\begin{theorem}\label{convexity_hm} For $n\ge 4$, $\delta\in (0,1)$, and $
\frac32<p\le 2$,
there exist $\epsilon_p=\epsilon(p,\delta)>0$ and $R_p=R(p,\delta)>0$ such that
if $u\in W^{2,2}(\Omega, N)$ is a biharmonic map satisfying either\\
(i) $\displaystyle\|\nabla ^2u\|_{M^{2,4}_{R_2}(\Omega)} \le\epsilon_2$, when $N$ is a
compact Riemannian manifold without boundary, or \\
(ii) $\displaystyle\|\nabla^2 u\|_{M^{p,2p}_{R_p}(\Omega)} \le\epsilon_p$,
when $N=\mathbb S^{L}$,\\
then
\begin{equation}\label{convex_ineq}
\int_{\Omega}|\Delta v|^2\ge \int_{\Omega}|\Delta u|^2+(1-\delta)\int_{\Omega}|\nabla^2(v-u)|^2
\end{equation}
holds for any $v\in W^{2,2}(\Omega, N)$ with $\displaystyle \left(v, \frac{\partial v}{\partial\nu}\right)
=\left(u,\frac{\partial u}{\partial\nu}\right)$ on $\partial\Omega$.
\end{theorem}

\pf First, it follows from Corollary 2.3 for $N=\mathbb S^L$ or Wang \cite{wang3}
that if $\epsilon_p>0$ is sufficiently small then $u\in C^\infty(\Omega, N)$, and
\begin{equation}\label{bi-blowup-rate}
\left|\nabla^m u(x)\right|\leq C\epsilon_p\left(\frac{1}{R_p^m}+\frac{1}{d^m(x,\partial\Omega)}\right),
\ \forall\ x\in \Omega, \ \forall\ m\geq 1.
\end{equation}
For $y\in N$, let $P^\perp(y):\mathbb R^{L+1}\to (T_yN)^\perp$ denote the orthogonal
projection from $\mathbb R^{L+1}$ to the normal space of $N$ at $y$. Since
$N$ is compact,  a simple geometric argument implies that there exists $C>0$ depending
on $N$ such that
\begin{equation}\label{orth_proj}
\Big|P^\perp(y)(z-y)\Big|\le C|z-y|^2, \ \forall z\in N.
\end{equation}
Since $$
\mathcal N_{\hbox{bh}}[u]\perp T_u N,$$
it follows from (\ref{orth_proj}) that multiplying (\ref{biharmonic-map}) by ($u-v$) and integrating over $\Omega$
yields
\begin{eqnarray}\label{unq-bi-map-eq4}
\int_{\Omega}\Delta u\cdot\Delta (u-v)
&=&\int_{\Omega}\mathcal N_{\hbox{bh}}[u]\cdot (u-v)\nonumber\\
&\lesssim& \int_{\Omega}[|\nabla u|^2|\nabla^2 u|+|\nabla^2 u|^2+|\nabla u||\nabla^3 u|]|u-v|^2\nonumber\\
&\lesssim& \epsilon_p^4 \int_\Omega \frac{|u-v|^2}{R_p^4}+\frac{|u-v|^2}{d^4(x,\partial\Omega)}\nonumber\\
&\lesssim& \epsilon_p\int_\Omega |\nabla^2(u-v)|^2,
\end{eqnarray}
where we choose $R_p\ge \epsilon_p$, apply (\ref{bi-blowup-rate}) and the Poincar\'e inequality and the Hardy
inequality (\ref{hardy-ineq}) during the last two steps.

It follows from (\ref{unq-bi-map-eq4}) that
\begin{equation}\label{unq-bi-map-eq8}
\int_{\Omega}|\Delta v|^2-\int_{\Omega}|\Delta u|^2-\int_{\Omega}|\Delta u-\Delta v|^2\\
=2\int_{\Omega}\Delta u\cdot\Delta (v-u)
\geq -C\epsilon_p\int_{\Omega}|\nabla^2(u-v)|^2.
\end{equation}
Since $(u-v)\in W^{2,2}_0(\Omega)$, we have that
$$\int_{\Omega}|\Delta u-\Delta v|^2=\int_{\Omega}|\nabla^2(u-v)|^2, $$
so that
$$\int_{\Omega}|\Delta v|^2-\int_{\Omega}|\Delta u|^2
\geq (1-C\epsilon_p)\int_{\Omega}|\nabla^2(u-v)|^2.$$
This yields (\ref{convex_ineq}), if $\epsilon_p>0$ is chosen so
that $C\epsilon_p\leq \delta$.
\endpf
\begin{corollary}\label{uniqueness} For
$n\ge 2$ and $\frac32<p\le 2$,  there exist $\epsilon_p>0$
and $R_p>0$ such that if $u_1, u_2\in W^{2,2}(\Omega, N)$ are biharmonic maps,
with $u_1-u_2\in W^{2,2}_0(\Omega,\mathbb R^{L+1})$,  satisfying either\\
(i) $\displaystyle\max_{i=1, 2}\|\nabla^2 u_i\|_{M^{2,4}_{R_2}(\Omega)}\le\epsilon_2$, when $N$ is a
compact Riemannian manifold without boundary, or \\
(ii) $\displaystyle\max_{i=1, 2}\|\nabla^2 u_i\|_{M^{p,2p}_{R_p}(\Omega)}
\le\epsilon_p$,
when $N=\mathbb S^{L}$, \\
then $u_1\equiv u_2$ in $\Omega$.
\end{corollary}

\pf   Choose $\delta=\frac12$, apply Theorem \ref{convexity_hm} to $u_1$ and $u_2$
by choosing sufficiently small $\epsilon_p>0$ and $R_p>0$. We have
$$
\int_{\Omega}|\Delta u_2|^2\ge \int_{\Omega}|\Delta u_1|^2+\frac12\int_{\Omega}|\nabla^2(u_2-u_1)|^2,
$$
and
$$
\int_{\Omega}|\Delta u_1|^2\ge \int_{\Omega}|\Delta u_2|^2+\frac12\int_{\Omega}|\nabla^2(u_1-u_2)|^2.
$$
Adding these two inequalities together yields $\displaystyle\int_{\Omega}|\nabla^2(u_1-u_2)|^2=0.$
This, combined with $u_1-u_2\in W^{2,2}_0(\Omega)$,
implies $u_1\equiv u_2$ in $\Omega$. \qed

\section {Uniqueness and convexity of heat flow of biharmonic maps}
\setcounter{equation}{0}
\setcounter{theorem}{0}
This section is devoted to the proof of uniqueness, convexity, and
 unique limit at $t=\infty$ for (\ref{biharmonic-flow}) of the heat flow of biharmonic maps, i.e.
Theorem \ref{biharmonic-flow-unique1},
Theorem \ref{biharmonic-flow-convex}, and Corollary \ref{biharmonic-flow-unique2}.

\medskip
\noindent{\bf Proof of Theorem \ref{biharmonic-flow-unique1}}.
First, by Theorem \ref{e-regularity},  we have that for $i=1, 2$,
$u_i\in C^{\infty}(\Omega\times(0,T),\mathbb S^L)$,
and
\begin{equation}\label{gradient_estimate10}
\Big|\nabla^m u_i(x,t)\Big|\leq C\epsilon_p
\left(\frac{1}{R_p^m}+\frac{1}{d^m(x,\partial\Omega)}
+\frac{1}{t^{\frac{m}{4}}}\right), \ \forall (x,t)\in\Omega\times (0, T), \
\forall\ m\ge 1.
\end{equation}
Set $w=u_1-u_2$. Then $w$ satisfies
\begin{equation}\label{bi-flow-uq-eq2}
\begin{cases}
\partial_tw+\Delta^2 w=\mathcal N_{\hbox{bh}}[u_1]-\mathcal N_{\hbox{bh}}[u_2]
&\ {\rm{in}}\ \Omega\times (0,T)\\
w=0 &\ {\rm{on}}\ \partial_p(\Omega\times (0,T))\\
\frac{\partial w}{\partial\nu}=0 &\ {\rm{on}}\ \partial \Omega\times (0,T).
\end{cases}
\end{equation}
Multiplying (\ref{bi-flow-uq-eq2})
by $w$ and integrating over $\Omega$,
by (\ref{orth_proj}),  (\ref{gradient_estimate10}),  the Poincar\'e inequality
and the Hardy inequality (\ref{hardy-ineq}),  we obtain that
\begin{eqnarray*}
\frac{d}{dt}\int_{\Omega}|w|^2+2\int_{\Omega}|\nabla^2 w|^2
&=&2\int_{\Omega}(\mathcal N_{\hbox{bh}}[u_1]-\mathcal N_{\hbox{bh}}[u_2])\cdot w\\
&\lesssim & \sum_{i=1}^2\int_{\Omega}
(|\nabla u_i|^2|\nabla^2 u_i|+|\nabla^2 u_i|^2+|\nabla u_i||\nabla^3 u_i|)|w|^2\\
&\lesssim& \epsilon_p^4\int_\Omega \frac{|w(x,t)|^2}{R_p^4}
+\frac{|w(x,t)|^2}{d^4(x,\partial\Omega)}+\frac{|w(x,t)|^2}{t}\\
&\lesssim& \epsilon_p\int_{\Omega}|\nabla^2w|^2
+\frac{\epsilon_p}{t}\int_{\Omega}|w|^2.
\end{eqnarray*}
If we choose $\epsilon_p>0$ sufficiently small and $R_p\ge \epsilon_p$,  then it holds
\begin{equation}\label{bi-flow-uq-eq3}
\frac{d}{dt}\int_{\Omega}|w|^2
\le\frac{C\epsilon_p}{t}\int_{\Omega}|w|^2.
\end{equation}
It follows from (\ref{bi-flow-uq-eq3}) that
\begin{eqnarray}\label{bi-flow-uq-eq4}
\frac{d}{dt}\Big(t^{-\frac12}\int_{\Omega}|w|^2\Big)
&=& t^{-\frac12}\frac{d}{dt}\int_{\Omega}|w|^2
-\frac12 t^{-\frac32}\int_{\Omega}|w|^2\nonumber\\
&\le& (C\epsilon-\frac12)t^{-\frac32}\int_{\Omega}|w|^2
\le 0.
\end{eqnarray}
Integrating this inequality from $0$ to $t$ yields
\begin{equation}\label{bi-flow-uq-eq5}
t^{-\frac12}\int_{\Omega}|w|^2
\le \lim_{t\downarrow 0^+}t^{-\frac12}\int_{\Omega}|w|^2.
\end{equation}
Since $w(\cdot,0)=0$, we have
$$w(x,t)=\int_0^t w_t(x,\tau)\,d\tau, \ {\rm{a.e.}}\  x\in \Omega,$$
so that, by the H\"older inequality,
$$t^{-\frac12}\int_{\Omega}|w(x,t)|^2
\le t^\frac12\int_0^t\int_{\Omega}|w_t|^2(x,\tau)\,dxd\tau
\le Ct^\frac12\rightarrow 0, \
{\rm{as}}\  t\downarrow 0^+.$$
This, combined with (\ref{bi-flow-uq-eq5}), implies
$w\equiv 0$ in $\Omega\times [0,T]$.
The proof is complete.
\endpf

\medskip
Now we want to prove Theorem \ref{biharmonic-flow-convex}
and Corollary \ref{biharmonic-flow-unique2}.
To do so, we need

\begin{lemma}\label{t-energy-mono}
Under the same assumptions as in Theorem \ref{biharmonic-flow-convex}, there exists
$T_0>0$ such that $\int_\Omega |\partial_t u(t)|^2$ is monotone decreasing for
$t\ge T_0$:
\begin{equation}\label{t-energy-mono1}
\int_\Omega |\partial_t u|^2(t_2)
+C\int_{\Omega\times [t_1, t_2]}|\nabla^2\partial_t u|^2
\le \int_\Omega |\partial_t u|^2(t_1), \ T_0\le t_1\le t_2\le T.
\end{equation}
\end{lemma}
\pf
For any sufficiently small $h>0$, set
$$u^{h}(x,t)=\frac{u(x,t+h)-u(x,t)}{h},\mbox{ }(x,t)\in \Omega\times(0,T-h).$$
Then $u^h\in L^2([0, T-h], W^{2,2}_0(\Omega))$,
$\partial_t u\in L^2(\Omega\times [0, T-h])$, and
$\displaystyle\lim\limits_{h\downarrow 0^+}\|u^h-\partial_t u\|_{L^{2}(\Omega\times[0,T-h])}=0.$
Since $u$ satisfies (\ref{biharmonic-flow}), we obtain
\begin{equation}\label{bi-flow-convex-eq3}
\partial_t u^h+\Delta^2 u^h=\frac{1}{h}\Big(\mathcal N_{\hbox{bh}}[u(t+h)]-\mathcal N_{\hbox{bh}}[u(t)]\Big).
\end{equation}
Multiplying (\ref{bi-flow-convex-eq3}) by $u^h$, integrating over $\Omega$, and applying (\ref{orth_proj})
and  (\ref{gradient_estimate10}), we have
\begin{eqnarray*}
\frac{d}{dt}\int_{\Omega}|u^h|^2+2\int_{\Omega}|\Delta u^h|^2
&\lesssim& \int_{\Omega}\left(|\mathcal N_{\hbox{bh}}[u(t+h)]|+|\mathcal N_{\hbox{bh}}[u(t)]|\right)|u^h|^2\\
&\lesssim& \int_{\Omega}\left(|\nabla^2u|^2+|\nabla u||\nabla^3 u|+|\nabla u|^2|\nabla^2 u||\right)
(t+h)|u^h|^2\\
&&+ \int_{\Omega}\left(|\nabla^2u|^2+|\nabla u||\nabla^3 u|+|\nabla u|^2|\nabla^2 u||\right)
(t)|u^h|^2\\
&\lesssim& \epsilon_p^4\int_\Omega \frac{|u^h|^2}{R_p^4}+\frac{|u^h|^2}{d^4(x,\partial\Omega)}
+\frac{|u^h|^2}{T_0}\\
&\lesssim& \epsilon_p \int_\Omega |\nabla^2 u^h|^2
\end{eqnarray*}
provided that we choose $R_p\ge \epsilon_p$ and $T_0\ge \epsilon_p$.
Since
$$\displaystyle\int_{\Omega}|\nabla^2u^h|^2=\int_{\Omega}|\Delta u^h|^2,$$
this implies
\begin{equation}\label{bi-flow-convex-eq4}
\frac{d}{dt}\int_{\Omega}|u^h|^2+2\int_{\Omega}|\nabla^2 u^h|^2
\leq \left(\frac12+C\epsilon_p\right)\int_\Omega |\nabla^2 u^h|^2.
\end{equation}
Choosing $\epsilon_p>0$ so that $C\epsilon_p\le \frac12$, integrating over $T_0\le t_1\le t_2\le T$, we have
\begin{equation}\label{bi-flow-convex-eq6}
\begin{split}
\int_{\Omega}|u^h|^2(t_2)+C\int_{t_1}^{t_2}\int_{\Omega}|\nabla^2 u^h|^2\leq \int_{\Omega}|u^h|^2(t_1).
\end{split}
\end{equation}
Sending $h\rightarrow0$, (\ref{bi-flow-convex-eq6})  yields (\ref{t-energy-mono1}). \qed

\medskip
Now we can show the monotonicity of  $E_2$-energy for heat flow of biharmonic maps
for $t\ge T_0$.

\begin{lemma}\label{2-energy-mono}
Under the same assumptions as in Theorem \ref{biharmonic-flow-convex}, there is
$T_0>0$ such that $\displaystyle\int_\Omega |\Delta u(t)|^2$ is monotone decreasing for
$t\ge T_0$:
\begin{equation}\label{2-energy-mono1}
\int_\Omega |\Delta u|^2(t_2)
+2\int_{\Omega\times [t_1, t_2]}|\partial_t u|^2
\le \int_\Omega |\Delta u|^2(t_1), \ T_0\le t_1\le t_2\le T.
\end{equation}
\end{lemma}
\pf For $\delta>0$, let $\eta_{\delta}\in C^{\infty}_0(\Omega)$ be such that
$$0\le \eta_\delta\le 1, \ \eta_{\delta}\equiv1 \ {\rm{for}}\ x\in\Omega\setminus
\Omega_{\delta},
\ {\rm{and}}\ |\nabla^m\eta_{\delta}|\leq C\delta^{-m},$$
where $\displaystyle\Omega_\delta=\{x\in\Omega: \ d(x,\partial\Omega)\le\delta\}$.
Multiplying (\ref{biharmonic-flow}) by $\partial_t u\eta_{\delta}^2$ and integrating
over $\Omega\times[t_1,t_2]$, we obtain
\begin{equation}\label{bi-flow-convex-eq7}
\begin{split}
&\int_{\Omega}|\Delta u(t_2)|^2\eta_{\delta}^2-\int_{\Omega}|\Delta u(t_1)|^2\eta_{\delta}^2+2\int_{t_1}^{t_2}\int_{\Omega}|\partial_tu|^2\eta_{\delta}^2\\
=&-4\int_{t_1}^{t_2}\int_{\Omega}\Delta u \cdot \partial_t u\left(|\nabla \eta_{\delta}|^2+\eta_{\delta}\Delta\eta_{\delta}\right)
-8\int_{t_1}^{t_2}\int_{\Omega}\Delta u\cdot\nabla \partial_t u\eta_{\delta}\nabla\eta_{\delta}.
\end{split}
\end{equation}
It suffices to show the right hand side of the above identity tends to $0$ as $\delta\rightarrow 0^+$.
By Lemma \ref{t-energy-mono}, we have that $\partial_t u \in L^2([T_0,T], W^{2,2}_0(\Omega))$
so that
\begin{equation}\label{bi-flow-convex-eq8}
\begin{split}
&\int_{t_1}^{t_2}\int_{\Omega}|\nabla\partial_t u|^2|\nabla\eta_{\delta}|^2
+|\partial_t u|^2\left(|\nabla\eta_{\delta}|^4+|\Delta\eta_{\delta}|^2\right)\\
\lesssim &\delta^{-2}\int_{t_1}^{t_2}\int_{\Omega_{\delta}}|\nabla \partial_t u|^2
+\delta^{-2}|\partial_t u|^2\\
\lesssim &\int_{t_1}^{t_2}\int_{\Omega_{\delta}}|\nabla^2 \partial_t u|^2\rightarrow 0,
\ {\rm{as}}\ \delta\rightarrow 0.
\end{split}
\end{equation}
This, combined with the H\"older inequality, implies that for $t_2\ge t_1\ge T_0$,
$$-4\int_{t_1}^{t_2}\int_{\Omega}\Delta u \cdot \partial_t u\left(|\nabla \eta_{\delta}|^2+\eta_{\delta}\Delta\eta_{\delta}\right)
-8\int_{t_1}^{t_2}\int_{\Omega}\Delta u\cdot\nabla \partial_t u\eta_{\delta}\nabla\eta_{\delta}\rightarrow 0,
\ {\rm{as}}\ \delta\rightarrow 0^+.$$
Thus (\ref{2-energy-mono1}) follows. \qed

\medskip
\noindent{\bf Proof of Theorem \ref{biharmonic-flow-convex}}.
First,  by Theorem \ref{e-regularity}, we have that
$u\in C^{\infty}(\Omega\times(0,T],\mathbb S^L)$, and
\begin{equation}\label{bi-flow-convex-eq1}
\Big|\nabla^mu(x,t)\Big|\leq C\epsilon_p
\left(\frac{1}{R_p^m}+\frac{1}{d^m(x,\partial\Omega)}+\frac{1}{t^{\frac{m}{4}}}\right),
\ \forall \ (x,t)\in\Omega\times (0, T), \ \forall\ m\ge 1.
\end{equation}
For $t_2>t_1\ge T_0$,  we have
\begin{equation}\label{bi-infty-eq2}
\begin{split}
&\int_{\Omega}|\Delta u(t_1)|^2-\int_{\Omega}|\Delta u(t_2)|^2-\int_{\Omega}|\Delta u(t_1)-\Delta u(t_2)|^2\\
=&2\int_{\Omega}\left(\Delta u(t_1)-\Delta u(t_2)\right)\Delta u(t_2)\\
=&-2\int_{\Omega}\left(u(t_1)-u(t_2)\right)u_t(t_2)\\
&+\int_{\Omega}\mathcal N_{\hbox{bh}}[u(t_2)]\cdot\left(u(t_1)-u(t_2)\right)\\
=&{I}+{II}.
\end{split}
\end{equation}
For $II$, applying (\ref{orth_proj}), we obtain
$$|\mathcal N_{\hbox{bh}}[u(t_2)]\cdot(u(t_1)-u(t_2))|\lesssim |\mathcal N_{\hbox{bh}}[u(t_2)]||u(t_1)-u(t_2)|^2.$$
Hence, by  (\ref{bi-flow-convex-eq1}), the Hardy inequality and the Poincar$\acute{\mbox{e}}$ inequality,
we have
\begin{equation}\label{bi-infty-eq3}
\begin{split}
|II|\lesssim &\epsilon_p^4\int_{\Omega}\left(\frac{1}{R_p^4}
+\frac{1}{d^4(x,\partial\Omega)}+\frac{1}{T_0}\right)\left|u(t_1)-u(t_2)\right|^2\\
\leq &C\epsilon_p\int_{\Omega}|\nabla^2(u(t_1)-u(t_2))|^2.
\end{split}
\end{equation}
For $I$,
by Lemma \ref{t-energy-mono}, we have
\begin{equation}\label{bi-infty-eq5}
\begin{split}
\Big\|\partial_t u(t_2)\Big\|_{L^2(\Omega)}^2
\leq \frac{1}{t_2-t_1}\int_{t_1}^{t_2}\int_{\Omega}|\partial_t u|^2.
\end{split}
\end{equation}
By the H$\ddot{\mbox{o}}$lder inequality and (\ref{2-energy-mono1}),  this implies
\begin{equation}\label{bi-infty-eq6}
\begin{split}
|I|\lesssim& \int_{\Omega}|\partial_t u(t_2)||u(t_1)-u(t_2)|\\
\lesssim& \left\|\partial_t u(t_2)\right\|_{L^2(\Omega)}
\left\|u(t_1)-u(t_2)\right\|_{L^2(\Omega)}\\
\leq& \sqrt{t_2-t_1}\left\|\partial_t u(t_2)\right\|_{L^2(\Omega)}
\left(\int_{\Omega\times [t_1, t_2]}|\partial_t u|^2\right)^{\frac{1}{2}}\\
\leq& \int_{\Omega\times [t_1, t_2]}|\partial_t u|^2
\leq\frac12\left[\int_{\Omega}|\Delta u(t_1)|^2-\int_{\Omega}|\Delta u(t_2)|^2\right].
\end{split}
\end{equation}
Putting (\ref{bi-infty-eq6}) and (\ref{bi-infty-eq3}) into (\ref{bi-infty-eq2}) implies (\ref{convexity2}).
This completes the proof.  \endpf\\

\noindent{\bf Proof of Corollary \ref{biharmonic-flow-unique2}}.
It follows from Lemma \ref{2-energy-mono} that $\displaystyle\int_\Omega |\Delta u(t)|^2$
is monotone decreasing for $t\geq T_0$. Hence
$$c=\lim\limits_{t\rightarrow\infty}\int_\Omega|\Delta u(t)|^2$$
exists and is finite.
Let $\{t_i\}$ be any increasing sequence such that $\lim\limits_{i\rightarrow\infty} t_i=+\infty$. Then
(\ref{convexity2}) implies that
$$\int_{\Omega}\Big|\nabla^2(u(t_{i+j})-u(t_i))\Big|^2
\leq C\Big[\int_{\Omega}|\Delta u(t_{i+j})|^2-\int_{\Omega}|\Delta u(t_i)|^2\Big]\rightarrow0,
\ {\rm{as}}\  i\rightarrow\infty, $$
for all $j\geq 1$. Thus there exists  $u_{\infty}\in W^{2,2}(\Omega,\mathbb S^L)$, with
$\displaystyle (u_{\infty}, \frac{\partial u_\infty}{\partial\nu})
=(u_0, \frac{\partial u_0}{\partial\nu})$ on $\partial \Omega$, such that
$$\lim\limits_{t\rightarrow\infty}\Big\|u(t)-u_{\infty}\Big\|_{W^{2,2}(\Omega)}=0.$$
Since (\ref{2-energy-mono1}) implies that there exists a sequence $t_i\rightarrow\infty$, such that
$$\lim\limits_{i\rightarrow\infty}\Big\|\partial_t u(t_i)\Big\|_{W^{2,2}(\Omega)}=0.$$
Thus $u_{\infty}\in W^{2,2}(\Omega, \mathbb S^L)$ is a biharmonic map.
For any $m\ge 1$, and any compact subset $K\subset\subset\Omega$, since
$$\Big\|u(t)\Big\|_{C^m(K)}\le C(n, m, K), \ \forall t\ge 1,$$
we conclude that
$$\lim_{t\rightarrow\infty}\Big\|u(t)-u_\infty\Big\|_{C^m(K)}=0,$$
and $u_\infty\in C^\infty(\Omega, \mathbb S^L)$.
This completes the proof.
\endpf

\section {Proof of Theorem \ref{biharmonic-flow-serrin}}
\setcounter{equation}{0}
\setcounter{theorem}{0}

In this section, we will prove Theorem \ref{biharmonic-flow-serrin} on both
smoothness and uniqueness for certain weak solutions of (\ref{biharmonic-flow}). First,
we would like to verify

\begin{proposition}\label{serrin_prop} For $n\ge 4$, $0<T<+\infty$,
suppose $u\in W^{1,2}_2(\Omega\times [0,T],N)$ is a weak solution of (\ref{biharmonic-flow}), with
the initial and boundary value $u_0\in W^{2,r}(\Omega, N)$ for some $\frac{n}2<r<+\infty$,
such that $\nabla^2 u\in L^q_tL^p_x(M\times [0,T])$ for some $p>\frac{n}2$ and $q<\infty$ satisfying (\ref{serrin_cond}). Then\\
(i) $\partial_t u\in L^{\frac{q}2}_tL^{\frac{p}2}_x(\Omega\times [0,T])$; and \\
(ii) for any $\epsilon>0$, there exists $R=R(u,\epsilon)>0$ such that
for any $1<s<\min\{\frac{p}2, \ \frac{q}2\}$,
\begin{equation}\label{small_norm6}
\sup\Big\{ r^{2s-(n+4)}\int_{P_r(x,t)\cap (\Omega\times [0,T])}(|\nabla^2 u|^s+r^{2s}|\partial_t u|^s)\ | \
(x,t)\in \Omega\times [0, T], \ 0<r\le R\Big\}\le \epsilon^s.
\end{equation}
\end{proposition}
\pf For simplicity, we will sketch the proof for $\Omega=\mathbb R^n$.
By the Duhamel formula, we have that $u(x,t)=u_1(x,t)+u_2(x,t)$, where
\begin{equation}\label{u1-eqn}
u_1(x,t)=\int_{\mathbb R^n} b(x-y, t) u_0(y),
\end{equation}
\begin{equation}\label{u2-eqn}
\begin{split}
&u_2(x,t)=\int_0^t\int_{\mathbb R^n} b(x-y, t-s)\mathcal N_{\hbox{bh}}[u](y,s)\\
&=\int_0^t\int_{\mathbb R^n} b(x-y, t-s)[\nabla\cdot(\nabla(A(u)(\nabla u,\nabla u))+2\Delta u\cdot\nabla(P(u)))
-\Delta u\cdot\Delta(P(u))](y,s).
\end{split}
\end{equation}
We proceed with two claims.\\
\noindent{\it Claim} 1. $\nabla^3 u\in L^{\frac{2q}3}_tL^{\frac{2p}3}_x(\mathbb R^n\times [0,T])$.
For $u_1$, we have
\begin{equation}\label{u3-eqn}
\nabla^3 u_1(x,t)=\int_{\mathbb R^n} \nabla_xb(x-y, t) \nabla^2 u_0(y).
\end{equation}
Direct calculations, using the property of the kernel function $b$, yield
\begin{equation}\label{u1-estimate}
\Big\|\nabla^3 u\Big\|_{L^{\frac{2q}3}_tL^{\frac{2p}3}_x(\mathbb R^n\times [0,T])}
\lesssim T^{\frac14(2-\frac{n}{r})}\Big\|\nabla^2 u_0\Big\|_{L^{r}(\mathbb R^n)}.
\end{equation}
For $u_2$, we have
\begin{eqnarray}\label{u4-eqn}
\nabla^3 u_2(x,t)&=&\int_0^t\int_{\mathbb R^n} \nabla_x^4b(x-y, t-s)\Big[\nabla(A(u)(\nabla u,\nabla u))+2\Delta u\cdot\nabla(P(u))\Big]\nonumber\\
&&-\int_0^t\int_{\mathbb R^n} \nabla_x^3b(x-y, t-s)\Delta u\cdot\Delta(P(u))(y,s)\nonumber\\
&=&M_1+M_2.
\end{eqnarray}
By the Nirenberg interpolation inequality, we have $\displaystyle\nabla u\in L^{2q}_tL^{2p}_x(\mathbb R^n\times [0,T])$. By the H\"older inequality,  we then have
$\displaystyle\nabla(A(u)(\nabla u,\nabla u))+2\Delta u\cdot\nabla(P(u)))\in L^{\frac{3q}2}_tL^{\frac{3p}2}_x(\mathbb R^n\times [0,T])$. Hence, by the Calderon-Zygmund $L^{\tilde q}_tL^{\tilde p}_x$-theory, we have
\begin{equation}\label{M1-estimate}
\begin{split}
\Big\|M_1\Big\|_{L^{\frac{2p}{3}}_{t}L^{\frac{2q}{3}}_x(\mathbb R^n\times [0,T])}
\lesssim & \Big\|\nabla(A(u)(\nabla u,\nabla u))+2\Delta u\cdot\nabla(P(u))\Big\|_{L^{\frac{2p}{3}}_{t}L^{\frac{2q}{3}}_x(\mathbb R^n\times [0,T])}\\
\lesssim &\Big\| \nabla u\Big\|_{L^{2p}_{t}L^{2q}_x(\mathbb R^n\times [0,T])}
\Big\| \nabla^2 u\Big\|_{L^{p}_{t}L^{q}_x(\mathbb R^n\times [0,T])}\\
\lesssim &1+\Big\| \nabla^2 u\Big\|_{L^{p}_{t}L^{q}_x(\mathbb R^n\times [0,T])}^2.
\end{split}
\end{equation}
For $M_2$, we have
$$|M_2|(x,t)\lesssim I_1\Big(|\nabla^2 u|^2+|\nabla u|^4\Big)(x,t),
 \ (x,t)\in \mathbb R^n\times [0,T].$$
Recall the following estimate of $I_1(\cdot)$ (see, for example, \cite{FJR} \S4):
\begin{equation}\label{I1-estimate}
\Big\|I_1(f)\Big\|_{L^{s_2}_tL^{r_2}_x(\mathbb R^n\times [0,T])}
\lesssim \Big\|f\Big\|_{L^{s_1}_tL^{r_1}_x(\mathbb R^n\times [0,T])},
\end{equation}
where $s_2\ge s_1$ and $r_2\ge r_1$ satisfy
\begin{equation}\label{index_cond}
\frac{n}{r_1}+\frac{4}{s_1}\le \frac{n}{r_2}+\frac{4}{s_2}+1.
\end{equation}
Applying (\ref{I1-estimate}) to $M_2$, we see that $M_2\in L^{\frac{2p}{3}}_{t}L^{\frac{2q}{3}}_x(\mathbb R^n\times [0,T])$,
and
\begin{equation}\label{M2-estimate}
\Big\|M_2\Big\|_{L^{\frac{2p}{3}}_{t}L^{\frac{2q}{3}}_x(\mathbb R^n\times [0,T])}
\lesssim 1+\Big\| \nabla^2 u\Big\|_{L^{p}_{t}L^{q}_x(\mathbb R^n\times [0,T])}^2.
\end{equation}
Combining these estimates of $\nabla^3 u_1, M_1,$ and $M_2$ yields {\it Claim} 1.

\medskip
\noindent{\it Claim} 2. $\nabla^4 u\in L^{\frac{q}2}_tL^{\frac{p}2}_x(\mathbb R^n\times [0,T])$.
It follows from {\it Claim} 1 that
$$\mathcal N_{\hbox{bh}}[u]=\displaystyle[\Delta(A(u)(\nabla u,\nabla u))+2\Delta u\cdot\nabla(P(u)))
-\Delta u\cdot\Delta(P(u))]\in L^{\frac{q}2}_tL^{\frac{p}2}_x(\mathbb R^n\times [0,T]).$$
Since
$$\nabla^4 u_2(x,t)=\int_0^t\int_{\mathbb R^n} \nabla^4_x b(x-y, t-s) \mathcal N_{\hbox{bh}}[u](y,s),$$
we can apply the Calderon-Zygmund $L^{\tilde q}_tL^{\tilde p}_x$-theory again to conclude that
$\nabla^4 u_2\in L^{\frac{q}2}_tL^{\frac{p}2}_x(\mathbb R^n\times [0,T])$.
For $u_1$, we have
$$\nabla^4u_1(x,t)=\int_{\mathbb R^n} \nabla_x^2b(x-y, t) \nabla^2 u_0(y).$$
Hence, by direct calculations, we have
$$
\Big\|\nabla^4 u_1\Big\|_{L^{\frac{q}2}_tL^{\frac{p}2}_x(\mathbb R^n\times [0,T])}
\lesssim T^{\frac14(2-\frac{n}r)}\Big\|\nabla^2 u_0\Big\|_{L^r(\mathbb R^n)}.
$$
Combining these two estimates yields {\it Claim} 2.

By (\ref{biharmonic-flow}), it is easy to see that $\partial_t u\in L^{\frac{q}2}_tL^{\frac{p}2}_x(\mathbb R^n\times [0,T])$.
In fact, we have
\begin{equation}\label{bi-flow-uq-eq15}
\begin{split}
\Big\|\partial_t u\Big\|_{L^{\frac{p}{2}}_{t}L^{\frac{q}{2}}_x(\mathbb R^n\times [0,T])}
\lesssim & \Big\|\mathcal N_{\hbox{bh}}[u]-\Delta^2 u\Big\|_{L^{\frac{p}{2}}_{t}L^{\frac{q}{2}}_x(\mathbb R^n\times [0,T])}\\
\lesssim &1+\Big\| \nabla^2 u\Big\|_{L^{p}_{t}L^{q}_x(\mathbb R^n\times [0,T])}^2+
T^{\frac14(2-\frac{n}r)}\Big\|\nabla^2 u_0\Big\|_{L^r(\mathbb R^n)}.
\end{split}
\end{equation}
This implies (i).

(ii) follows from (i) and the H\"older inequality. In fact,
for any $1<s<\min\{\frac{p}{2}, \frac{q}2\}$, it holds
$$\Big(r^{2s-(n+4)}\int_{P_r(x,t)\cap (\Omega\times [0, T])} |\nabla ^2u|^s\Big)^{\frac1{s}}
\le \Big\|\nabla^2 u\Big\|_{L^{q}_tL^{p}_x(P_r(x,t)\cap (\Omega\times [0, T]))},$$
and
$$\Big(r^{4s-(n+4)}\int_{P_r(x,t)\cap (M\times [0, T])} |\partial_t u|^s\Big)^{\frac1{s}}
\le \Big\|\partial_t u\Big\|_{L^{\frac{q}2}_tL^{\frac{p}2}_x(P_r(x,t)\cap (\Omega\times [0, T]))}.$$
These two inequalities clearly imply (\ref{small_norm6}), provided that $R=R(u,\epsilon)>0$
is chosen sufficiently small.
\qed

\medskip
Now we prove an $\epsilon$-regularity property for certain solutions of (\ref{biharmonic-flow}).

\begin{theorem} \label{e-regularity-serrin} There exists $\epsilon_0>0$ such that
if $u\in W^{1,2}_2(P_1, N)$, with $\nabla^2 u\in L^q_tL^p_x(P_1)$ for some
$q\ge \frac{n}2$ and $p\le \infty$ satisfying (\ref{serrin_cond}), is a weak solution
of (\ref{biharmonic-flow}) and satisfies
\begin{equation}
\label{e-cond-serrin}
\Big\|\nabla^2 u\Big\|_{L^{q}_tL^p_x(P_1)}\leq \epsilon_0,
\end{equation}
then $u\in C^{\infty}(P_{\frac{1}{2}},N)$ and
\begin{equation}\label{grad_estimate10}
\|\nabla^m u\|_{C^0(P_\frac12)}\le C(m, p, q, n)\|\nabla^2 u\|_{L^q_tL^p_x(P_1)}, \ \forall\ m\ge 1.
\end{equation}
\end{theorem}

Before proving this theorem, we recall the Serrin type inequalities (see \cite{Serrin}) and
Adams' estimates of Riesz potential between Morrey spaces in $(\mathbb R^{n+1},\delta)$.

\begin{lemma}\label{bi-serrin-lemma}{\it Assume $p\ge \frac{n}2$ and $q\le\infty$ satisfy (\ref{serrin_cond}).
For any $f\in L^q_tL^p_x(\Omega\times [0,T])$, $g\in L^2_tW^{2,2}_x(\Omega\times[0,T])$,
and $h\in L^2_tW^{1,2}_x(\Omega\times [0,T])$, we have
\begin{equation}\label{bi-gen-serrin-eq1}
\int_{\Omega\times [0,T]}|f||g||h|\lesssim\|h\|_{L^2(\Omega\times[0,T])}
\|g\|^{\frac{n}{2p}}_{L^2_tW^{2,2}_x(\Omega\times[0,T])}
\left(\int_0^T\|f\|_{L^p(\Omega)}^q\|g\|_{L^2(\Omega)}^2\right)^{\frac{1}{q}},
\end{equation}
and
\begin{equation}\label{bi-gen-serrin-eq2}
\int_{\Omega\times [0,T]}|f||\nabla g||h|\lesssim\|h\|_{L^2_tW^{1,2}_x(\Omega\times[0,T])}
\|g\|^{\frac{n}{2p}}_{L^2_tW^{2,2}_x(\Omega\times[0,T])}
\left(\int_0^T\|f\|_{L^p(\Omega)}^q\|g\|_{L^2(\Omega)}^2\right)^{\frac{1}{q}}.
\end{equation}
}
\end{lemma}

\pf For convenience, we sketch the proof here. By the H$\ddot{\mbox{o}}$lder inequality, we have
\begin{equation}\label{bi-gen-serrin-eq3}
\int_{\Omega}|f||g||h|\leq \|f\|_{L^p(\Omega)}\|g\|_{L^r(\Omega)}\|h\|_{L^2(\Omega)},
\end{equation}
where $\displaystyle \frac{1}{p}+\frac{1}{r}=\frac{1}{2}.$
It follows from (\ref{serrin_cond}) that $2\leq r\leq \frac{2n}{n-4}$.
Hence by the Sobolev inequality we have
\begin{equation}\label{bi-gen-serrin-eq4}
\|g\|_{L^r(\Omega)}\leq\|g\|^{\frac{2}{q}}_{L^2(\Omega)}
\|g\|_{L^{\frac{2n}{n-4}}(\Omega)}^{\frac{2n}{p}}
\lesssim \|g\|^{\frac{2}{q}}_{L^2(\Omega)}\|g\|^{\frac{n}{2p}}_{W^{2,2}(\Omega)}.
\end{equation}
Putting (\ref{bi-gen-serrin-eq4}) into (\ref{bi-gen-serrin-eq3}) yields
\begin{equation}\label{bi-gen-serrin-eq5}
\int_{\Omega}|f||g||h|\lesssim\|f\|_{L^p(\Omega)}
\|g\|^{\frac{2}{q}}_{L^2(\Omega)}\|g\|^{\frac{n}{2p}}_{W^{2,2}(\Omega)}\|h\|_{L^2(\Omega)}.
\end{equation}
Since $\displaystyle\frac{1}{q}+\frac{n}{4p}+\frac{1}{2}=1,$
(\ref{bi-gen-serrin-eq1}) follows by integrating over $[0,T]$ and the H$\ddot{\mbox{o}}$lder inequality.

To see (\ref{bi-gen-serrin-eq2}), note that the H$\ddot{\mbox{o}}$lder inequality implies
\begin{equation}\label{bi-gen-serrin-eq6}
\int_{\Omega}|f||\nabla g||h|\leq \|f\|_{L^p(\Omega)}\|\nabla g\|_{L^s(\Omega)}\|h\|_{L^{\frac{2n}{n-2}}(\Omega)}
\end{equation}
where $\displaystyle\frac{1}{p}+\frac{1}{s}+\frac{n-2}{2n}=1.$
Since 
$$\displaystyle \frac{1}{s}=\frac{1}{n}+\frac{n}{2p}\left(\frac{1}{2}-\frac{2}{n}\right)
+ \left(1-\frac{n}{2p}\right)\frac{1}{2},$$ 
the Nirenberg interpolation inequality implies
\begin{equation}\label{bi-gen-serrin-eq7}
\|\nabla g\|_{L^s(\Omega)}\lesssim\|g\|^{\frac{2}{q}}_{L^2(\Omega)}\|g\|^{\frac{n}{2p}}_{W^{2,2}(\Omega)}.
\end{equation}
Putting (\ref{bi-gen-serrin-eq7}) into (\ref{bi-gen-serrin-eq6}) and using the Sobolev inequality, we obtain
\begin{equation}\label{bi-gen-serrin-eq8}
\int_{\Omega}|f||\nabla g||h|\lesssim\|f\|_{L^p(\Omega)}\|g\|^{\frac{2}{q}}_{L^2(\Omega)}
\|g\|^{\frac{n}{2p}}_{W^{2,2}(\Omega)}\|h\|_{W^{1,2}(\Omega)}.
\end{equation}
Since $\displaystyle\frac{1}{q}+\frac{n}{4p}+\frac{1}{2}=1,$ (\ref{bi-gen-serrin-eq2}) follows
by integration on $[0,T]$ and the H$\ddot{\mbox{o}}$lder inequality. \endpf\\

Now we state Adams' estimate for the Riesz potentials  on $(\mathbb R^{n+1}, \delta)$.
Since its proof is exactly  the same argument as in Huang-Wang (\cite{huang-wang1} Theorem 3.1),
 we  skip it here.

\begin{proposition}\label{bi-riesz-theorem}{\it
 (i) For any $\beta>0$, $0 <\lambda\leq n+4$, $1<p<\frac{\lambda}{\beta}$, if $f\in L^p(\mathbb R^{n+1})\cap M^{p,\lambda}(\mathbb R^{n+1})$,
then $I_{\beta}(f)\in L^{\tilde{p}}(\mathbb R^{n+1})\cap M^{\tilde{p},\lambda}(\mathbb R^{n+1})$, where $\tilde{p}=\frac{p\lambda}{\lambda-p\beta}$. Moreover,
\begin{equation}
\label{reg-eqn3}
\|I_{\beta}(f)\|_{L^{\tilde{p}}(\mathbb R^{n+1})}\leq C\|f\|^{\frac{\beta p}{\lambda}}_{M^{p,\lambda}(\mathbb R^{n+1})}\|f\|^{1-\frac{\beta p}{\lambda}}_{L^{p}(\mathbb R^{n+1})}
\end{equation}
\begin{equation}
\label{reg-eqn4}
\|I_{\beta}(f)\|_{M^{\tilde{p},\lambda}(\mathbb R^{n+1})}\leq C\|f\|_{M^{p,\lambda}(\mathbb R^{n+1})}.
\end{equation}\\
(ii) For any $0<\beta<\lambda\leq n+4$, if $f\in L^1(\mathbb R^{n+1})\cap M^{1,\lambda}(\mathbb R^{n+1})$,
then $f\in L^{\frac{\lambda}{\lambda-\beta},*}(\mathbb R^{n+1})\cap M_*^{\frac{\lambda}{\lambda-\beta},\lambda}(\mathbb R^{n+1})$. Moreover,
\begin{equation}
\label{reg-eqn5}
\|I_{\beta}(f)\|_{L^{\frac{\lambda}{\lambda-\beta},*}(\mathbb R^{n+1})}\leq C\|f\|^{\frac{\beta}{\lambda}}_{M^{1,\lambda}(\mathbb R^{n+1})}\|f\|^{1-\frac{\beta}{\lambda}}_{L^{1}(\mathbb R^{n+1})}
\end{equation}
\begin{equation}
\label{reg-eqn6}
\|I_{\beta}(f)\|_{M_*^{\frac{\lambda}{\lambda-\beta},\lambda}(\mathbb R^{n+1})}\leq C\|f\|_{M^{1,\lambda}(\mathbb R^{n+1})}.
\end{equation}
}
\end{proposition}

\bigskip
\noindent{\bf Proof of Theorem \ref{e-regularity-serrin}}. The proof  is based on three claims.

\medskip
\noindent{\it Claim 1}.   For any $0<\alpha<1$,  we have that $\displaystyle\nabla^2 u\in M^{2,4-4\alpha}(P_\frac34)$,
and
\begin{equation}\label{bi-gen-serrin-eq19}
\Big\|\nabla^2 u\Big\|_{M^{2,4-4\alpha}(P_\frac34)}
\le C \Big\|\nabla^2 u\Big\|_{L^q_tL^p_x(P_1)}.
\end{equation}
For any $0<r\leq \frac{1}{4}$ and $z_0=(x_0,t_0)\in P_{\frac{3}{4}}$, by (\ref{e-cond-serrin}) we have
\begin{equation}
\label{bi-gen-serrin-eq9}
\|\nabla^2 u\|_{L^{q}_tL^p_x(P_r(z_0))}\leq \epsilon.
\end{equation}
Let $v:P_r(z_0)\rightarrow \mathbb R^{L+1}$ solve
\begin{equation}
\label{bi-gen-serrin-eq10}
\left\{
\begin{split}
\partial_t v+\Delta^2 v=&0\quad\mbox{in }P_r(z_0)\\
 v=&u\quad\mbox{on }\partial_pP_r(z_0)\\
 \frac{\partial v}{\partial\nu}=& \frac{\partial u}{\partial\nu}\ \mbox{on }\partial B_r(x_0)\times(t_0-r^4, t_0].
\end{split}
\right.
\end{equation}
Set $w=u-v$. Multiplying (\ref{bi-gen-serrin-eq10}) and (\ref{biharmonic-flow}) by $w$, subtracting the resulting equations and integrating over $P_r(z_0)$, we obtain
\begin{equation}
\label{bi-gen-serrin-eq11}
\begin{split}
&\sup\limits_{t_0-r^4\leq t\leq t_0}\int_{B_r(x_0)}|w|^2(t)+2\int_{P_r(z_0)}|\nabla^2w|^2\\
=& |\int_{P_r(z_0)}\mathcal N_{\hbox{bh}}[u]\cdot w|\\
=& |\int_{P_r(z_0)}-\nabla(A(u)(\nabla u,\nabla u))\nabla w
-\left<\Delta u,\Delta(P(u))\right>w
-2\left<\Delta u,\nabla(P(u))\right>\nabla w|\\
\lesssim& \int_{P_r(z_0)}|\nabla^2u|^2|w|+\int_{P_r(z_0)}|\nabla u||\nabla^2u||\nabla w|\\
=&I+II.
\end{split}
\end{equation}
For $I$, we can apply (\ref{bi-gen-serrin-eq1}) to get
\begin{equation}
\label{bi-gen-serrin-eq12}
\begin{split}
|I|\lesssim\|\nabla^2 u\|_{L^2(P_r(z_0))}
\|w\|^{\frac{n}{2p}}_{L^2_tW^{2,2}_x(P_r(z_0))}
\left(\int_{t_0-r^4}^{t_0}\|\nabla^2 u\|_{L^p(B_r(x_0))}^q\|w\|_{L^2(B_r(x_0))}^2\right)^{\frac{1}{q}}.
\end{split}
\end{equation}
For $II$,  by (\ref{bi-gen-serrin-eq2}), we have
\begin{equation}
\label{bi-gen-serrin-eq13}
\begin{split}
|II|\lesssim\|\nabla u\|_{L^2_tW^{1,2}_x(P_r(z_0)}
\|w\|^{\frac{n}{2p}}_{L^2_tW^{2,2}_x(P_r(z_0))}
\left(\int_{t_0-r^4}^{t_0}\|\nabla^2 u\|_{L^p(B_r(x))}^q\|w\|_{L^2(B_r(x_0))}^2\right)^{\frac{1}{q}}.
\end{split}
\end{equation}
Putting (\ref{bi-gen-serrin-eq12}) and (\ref{bi-gen-serrin-eq13}) into (\ref{bi-gen-serrin-eq11}) and applying the Poincar$\acute{\mbox{e}}$ inequality, we obtain
\begin{equation}
\label{bi-gen-serrin-eq14}
\begin{split}
&\sup\limits_{t_0-r^4\leq t\leq t_0}\int_{B_r(x_0)}|w|^2(t)+2\int_{P_r(z_0)}|\nabla^2 w|^2\\
\lesssim&\begin{cases} \|\nabla u\|_{L^2_tW^{1,2}_x(P_r(z_0))}
\|\nabla ^2w\|^{\frac{n}{2p}}_{L^2(P_r(z_0))}
\left(\int_{t_0-r^4}^{t_0}\|\nabla^2 u\|_{L^p(B_r(x_0))}^q\|w\|_{L^2(B_r(x_0))}^2\right)^{\frac{1}{q}},\ & q<\infty,\\
\|\nabla u\|_{L^2_tW^{1,2}_x(P_r(z_0))}
\|\nabla ^2w\|_{L^2(P_r(z_0))}\|\nabla^2 u\|_{L^\infty_t L^{\frac{n}2}_x(B_r(x_0))}, \ & q=\infty.
\end{cases}
\end{split}
\end{equation}
Since $\displaystyle\|\nabla^2 u\|_{L^q_tL^p_x(P_r(z_0))}\le\epsilon$,  we obtain, by the Young inequality,
\begin{eqnarray}
&\sup\limits_{t_0-r^4\leq t\leq t_0}\int_{B_r(x_0)}|w|^2(t)+2\int_{P_r(z_0)}|\nabla^2 w|^2\nonumber\\
\leq&\begin{cases}
\|\nabla^2 w\|^2_{L^2(P_r(z_0))}+\epsilon\|\nabla u\|^2_{L^2_tW^{1,2}_x(P_r(z_0))}
+C\epsilon^{\frac{p}2}\sup\limits_{t_0-r^4\leq t\leq t_0}\|w\|_{L^2(B_r(x_0))}^2, & q<\infty,\\
\|\nabla^2 w\|_{L^2(P_r(z_0))}^2+
C\|\nabla^2 u\|^2_{L^\infty_t L^{\frac{n}2}_x(B_r(x_0))}\|\nabla u\|_{L^2_tW^{1,2}_x(P_r(z_0))}^2, \ & q=\infty.
\end{cases}
\end{eqnarray}
By choosing $\epsilon>0$ sufficiently small, this implies
\begin{equation}
\label{bi-gen-serrin-eq16}
\int_{P_r(z_0)}|\nabla^2 w|^2
\lesssim \epsilon\int_{P_r(z_0)}|\nabla u|^2+|\nabla^2 u|^2.
\end{equation}
Since $N$ is compact and $u$ maps into $N$, $|u|\le C_N$.
Hence, by the Nirenberg interpolation inequality, we have
\begin{equation}
\label{bi-gen-serrin-eq17}
\int_{P_r(z_0)}|\nabla u|^2\lesssim\int_{P_r(z_0)}|\nabla^2u|^2+r^{n+4}.
\end{equation}
Combining (\ref{bi-gen-serrin-eq17}) with (\ref{bi-gen-serrin-eq16}), we have
\begin{equation}
\label{bi-gen-serrin-eq18}
\int_{P_r(z_0)}|\nabla^2 w|^2
\lesssim\epsilon\int_{P_r(z_0)}|\nabla^2u|^2+\epsilon r^{n+4}.
\end{equation}
By the standard estimate on $v$, we have
\begin{equation}\label{v-estimate}
(\theta r)^{-n}\int_{P_{\theta r}(z_0)}|\nabla^2 v|^2
\lesssim\theta^4r^{-n}\int_{P_r(z_0)}|\nabla^2v|^2, \ \forall\ \theta\in(0,1).
\end{equation}
Combining (\ref{bi-gen-serrin-eq18}) with (\ref{v-estimate}), we obtain
\begin{equation}
\label{bi-gen-serrin-eq20}
(\theta r)^{-n}\int_{P_{\theta r}(z_0)}|\nabla^2 u|^2
\leq C\left(\theta^4+\theta^{-n}\epsilon\right)r^{-n}\int_{P_r(z_0}|\nabla^2 u|^2+C\epsilon \theta^{-n}r^{4},
\ \forall\ \theta\in(0,1).
\end{equation}
For any $0<\alpha<1$, choose $0<\theta<1$ and $\epsilon$ such that
$$C\theta^4\leq \frac{1}{2}\theta^{4\alpha}\ {\rm{and}}\
\epsilon\leq \min\left\{\left(\frac{1}{2C}\right)^{\frac{2}{p}}, \frac{\theta^{4\alpha+n}}{2C}\right\}.$$
Therefore,  for any $(z_0)\in P_{\frac{3}{4}}$ and $0<r\leq \frac{1}{4}$, it holds
\begin{equation}
\label{bi-gen-serrin-eq21}
(\theta r)^{-n}\int_{P_{\theta r}(x,t)}|\nabla^2 u|^2
\leq \theta^{4\alpha}r^{-n}\int_{P_r(x,t)}|\nabla^2 u|^2+\theta^{4\alpha}r^{4}.
\end{equation}
It is  standard that iterating (\ref{bi-gen-serrin-eq21}) implies
\begin{equation}
\label{bi-gen-serrin-eq23}
r^{-n}\int_{P_{ r}(z_0)}|\nabla^2 u|^2
\leq Cr^{4\alpha}\left(\int_{P_1}|\nabla^2 u|^2+1\right)
\end{equation}
for any $z_0\in P_{\frac{3}{4}}$ and $0<r\leq\frac{1}{4}$.
(\ref{bi-gen-serrin-eq23}) implies
that $\nabla^2u \in M^{2, 4-4\alpha}(P_\frac34)$,
and the estimate (\ref{bi-gen-serrin-eq19}) holds. This proves {\it Claim} 1.

\medskip
\noindent{\it Claim 2}. For any $1<\beta<+\infty$, $\nabla^2 u\in L^\beta(P_\frac9{16})$, and
\begin{equation}
\Big\|\nabla^2 u\Big\|_{L^\beta(P_\frac9{16})}
\lesssim \Big\|\nabla^2 u\Big\|_{L^q_tL^p_x(P_1)}^2.\label{morrey6}
\end{equation}
This can be proven by estimates of Riesz potentials between Morrey spaces.
To do so, let $\eta \in C^{\infty}_0(P_1)$ be such that
$$0\leq \eta\leq 1, \ \eta\equiv1 \ {\rm{in}}\ P_{\frac{5}{8}},
\ |\eta_t|+\sum\limits_{m=1}^4|\nabla^m \eta|\leq C.$$
Let $Q:\mathbb R^n\times [-1, \infty]\to\mathbb R^{L+1}$ solve
\begin{eqnarray}\label{Q-eqn}
\partial_t Q+\Delta^2 Q &=& \nabla\cdot\Big(\eta^2 \nabla(A(u)(\nabla u, \nabla u))
+2\eta^2\langle\Delta u, \nabla(P(u))\rangle\Big)-\eta^2\langle\Delta u, \Delta(P(u))\rangle\\
Q\Big|_{t=-1}&=& 0. \nonumber
\end{eqnarray}
Set
$$J_1= \nabla\cdot\Big(\eta^2 \nabla(A(u)(\nabla u, \nabla u))
+2\eta^2\langle\Delta u, \nabla(P(u))\rangle\Big)\ \ {\rm{and}}\
\  J_2= -\eta^2\langle\Delta u, \Delta(P(u))\rangle.$$
By  the Duhamel formula, we have, for $(x,t)\in \mathbb R^n\times (-1, \infty)$,
\begin{equation}
\label{bi-gen-serrin-eq26}
\begin{split}
\nabla^2Q(x,t)=&\int_{\mathbb R^n\times [-1, t]}\nabla^2_xb(x-y,t-s)\left(J_1+J_2\right)(y,s)\\
=&\int_{\mathbb R^{n}\times [-1, t]}\nabla_x^3b(x-y,t-s)\Big(\eta^2 \nabla(A(u)(\nabla u, \nabla u))
+2\eta^2\langle\Delta u, \nabla(P(u))\rangle\Big)(y,s)\\
&-\int_{\mathbb R^n\times [-1, t]}\nabla^2_xb(x-y, t-s)\eta^2\langle\Delta u, \Delta(P(u))\rangle(y,s)\\
=&K_1(x,t)+K_2(x,t).
\end{split}
\end{equation}
It is clear that for $(x,t)\in \mathbb R^n\times (-1, \infty)$,
$$|K_1|(x,t)\lesssim I_1\Big(\eta^2 (|\nabla u|^3+|\nabla u||\nabla^2 u|)\Big)(x,t),
\ |K_2|(x,t)\le I_2\Big(\eta^2 (|\nabla^2 u|^2+|\nabla u|^4)\Big)(x,t).$$
It follows from (\ref{bi-gen-serrin-eq19}) and the Nirenberg interpolation inequality
that $\nabla u\in M^{4,4-4\alpha}(P_\frac34)$ and
\begin{equation}\label{morrey2}
\Big\|\nabla u\Big\|_{M^{4,4-4\alpha}(P_\frac34)}\lesssim \Big\|\nabla^2 u\Big\|_{L^q_tL^p_x(P_1)}.
\end{equation}
Hence, by the H\"older inequality, we have that for any $0<\alpha_1, \alpha_2<1$,
$$\displaystyle \eta^2 (|\nabla u|^3+|\nabla u||\nabla^2 u|)\in M^{\frac43, 4-4\alpha_1}(\mathbb R^{n+1})
\ {\rm{and}}\ \displaystyle \eta^2 (|\nabla^2 u|^2+|\nabla u|^4)\in M^{1, 4-4\alpha_2}(\mathbb R^{n+1}),$$
and
\begin{eqnarray}\label{morrey3}
\Big\|\eta^2 (|\nabla u|^3+|\nabla u||\nabla^2 u|)\Big\|_{M^{\frac43, 4-4\alpha_1}(\mathbb R^{n+1})}
&\lesssim& \Big\|\nabla u\Big\|_{M^{4,4-4\alpha_1}(P_\frac34)}\Big\|\nabla^2 u\Big\|_{M^{2,4-4\alpha_1}(P_\frac34)}
\nonumber\\
&\lesssim& \Big\|\nabla^2 u\Big\|_{L^q_tL^p_x(P_1)}^2,
\end{eqnarray}
\begin{eqnarray}\label{morrey4}
\Big\|\eta^2 (|\nabla^2 u|^2+|\nabla u|^4)\Big\|_{M^{1, 4-4\alpha_2}(\mathbb R^{n+1})}
&\lesssim& \Big\|\nabla u\Big\|_{M^{4,4-4\alpha_2}(P_\frac34)}+\Big\|\nabla^2 u\Big\|_{M^{2,4-4\alpha_2}(P_\frac34)}
\nonumber\\
&\lesssim& \Big\|\nabla^2 u\Big\|_{L^q_tL^p_x(P_1)}^2.
\end{eqnarray}
Now applying Proposition \ref{bi-riesz-theorem}, we conclude that
$$\displaystyle K_1\in M^{\frac{4-4\alpha_1}{2-3\alpha_1},4-4\alpha_1}\cap  L^{\frac{4-4\alpha_1}{2-3\alpha_1}}(\mathbb R^{n+1}),
\ \ \displaystyle K_2\in M_*^{\frac{2-2\alpha_2}{1-2\alpha_2},4-4\alpha_2}\cap
L^{\frac{2-2\alpha_2}{1-2\alpha_2}, *}(\mathbb R^{n+1}),$$
and
\begin{equation}\label{morrey5}
\Big\|K_1\Big\|_ {M^{\frac{4-4\alpha_1}{2-3\alpha_1},4-4\alpha_1}(\mathbb R^{n+1})}
+\Big\|K_2\Big\|_{M_*^{\frac{2-2\alpha_2}{1-2\alpha_2},4-4\alpha_2}(\mathbb R^{n+1})}
\lesssim  \Big\|\nabla^2 u\Big\|_{L^q_tL^p_x(P_1)}^2.
\end{equation}
Sending $\alpha_1\uparrow \frac{2}{3}$ and $\alpha_2\uparrow \frac12$,
we obtain that for any $1<\beta<+\infty$,
$K_1, K_2\in L^\beta(\mathbb R^{n+1})$, and
\begin{equation}\|K_1\|_{L^\beta(\mathbb R^{n+1})}
+\|K_2\|_{L^\beta(\mathbb R^{n+1})}\lesssim \Big\|\nabla^2 u\Big\|_{L^q_tL^p_x(P_1)}^2.\label{morrey6}
\end{equation}
This implies that for any $1<\beta<+\infty$, $\nabla^2 Q\in L^\beta(\mathbb R^{n+1})$, and
\begin{equation}\Big\|\nabla^2 Q\Big\|_{L^\beta(\mathbb R^{n+1})}
\lesssim \Big\|\nabla^2 u\Big\|_{L^q_tL^p_x(P_1)}^2.\label{morrey6}
\end{equation}
Since ($u-Q$) solves
$$\Big(\partial_t+\Delta^2\Big)(u-Q)=0 \ {\rm{in}}\ P_\frac58,$$
it follows that for any $1<\beta<+\infty$, $\nabla^2 u\in L^\beta(P_\frac9{16})$, and
\begin{equation}
\Big\|\nabla^2 u\Big\|_{L^\beta(P_\frac9{16})}
\lesssim \Big\|\nabla^2 u\Big\|_{L^q_tL^p_x(P_1)}^2.\label{morrey7}
\end{equation}
This implies (\ref{morrey6}). Hence {\it Claim} 2 is proven.

\medskip
\noindent{\it Claim 3}. $u\in C^\infty(P_\frac12, N)$ and (\ref{grad_estimate10}) holds.
It follows from (\ref{morrey6}) that for any $1<\beta<+\infty$, there exist $f, g\in L^\beta(P_{\frac9{16}})$ such that
(\ref{biharmonic-flow}) can be written as
$$\displaystyle(\partial_t +\Delta^2 )u=\nabla\cdot f+g.$$
Thus, by the $L^p$-theory of higher-order parabolic equations, we conclude that
$\nabla^3 u\in L^\beta(P_\frac{17}{32})$. Applying the $L^p$-theory again, we would obtain that
$\partial_t u,\nabla^4 u\in L^\beta(P_{\frac{33}{64}})$.  Taking derivatives of the equation
(\ref{biharmonic-flow}) and repeating this argument, we can conclude
that $u\in C^\infty(P_\frac12, N)$, and the estimate (\ref{grad_estimate10}) holds.
Putting together these three claims completes the proof.
\endpf

\medskip
\noindent{\bf Proof of Theorem \ref{biharmonic-flow-serrin}}. Let $\epsilon_0>0$ be given by Theorem \ref{e-regularity-serrin}. Since $p>\frac{n}2$ and $q<\infty$, there exists $T_0>0$ such that
\begin{equation}\label{serrin_small1}
\max_{i=1,2}\|\nabla^2 u_i\|_{L^q_tL^p_x(\Omega\times [0,T_0])}\le\epsilon_0.
\end{equation}
This implies that for any $x_0\in\Omega$ and $0<t_0\le T_0$, if $R_0=\min\{d(x_0,\partial\Omega), t_0^{\frac14}\}>0$,
then
\begin{equation}\label{serrin_small2}
\max_{i=1,2}\|\nabla^2 u_i\|_{L^q_tL^p_x(P_{R_0}(z_0))}\le\epsilon_0.
\end{equation}
Hence by suitable scalings of the estimate of Theorem \ref{e-regularity-serrin}, we have that for $i=1,2$, $u_i\in C^\infty(P_{\frac{R_0}2}(z_0), N)$
and
\begin{equation}\label{gradient_estimate60}
\Big|\nabla^m u_i\Big|(x_0,t_0)\lesssim \epsilon_0 \left(\frac{1}{d^m(x_0,\partial\Omega)}+\frac{1}{t_0^{\frac{m}4}}\right).
\end{equation}
Using (\ref{gradient_estimate60}), the same proof of Theorem \ref{biharmonic-flow-unique1} implies
that $u_1\equiv u_2$ in $\Omega\times [0,T_0]$. Repeating this argument on the interval $[T_0, T]$
yields  $u_1\equiv u_2$ in $\Omega\times [0,T]$. \qed

\medskip
\noindent{\bf Proof of Corollary \ref{4d_unique}}.  Let $\epsilon_0>0$ be given by Theorem \ref{e-regularity-serrin}.
Since $u_0\in W^{2,2}(\Omega, N)$, by the absolute continuity of $\displaystyle\int |\nabla^2 u_0|^2$
there exists $r_0>0$ such that
\begin{equation}\label{serrin_small3}
\max_{x\in\Omega}\int_{B_{r_0}(x)\cap\Omega}|\nabla^2 u_0|^2\le\frac{\epsilon_0^2}2.
\end{equation}
Choosing $\epsilon_1\le \frac{\epsilon_0^2}2$ and applying (\ref{no_jump}), we conclude that
there exists $0<t_0\le r_0^4$ such that
\begin{equation}\label{serrin_small4}
\max_{x\in\Omega, 0\le t\le t_0}\int_{B_{r_0}(x)\cap\Omega}|\nabla^2 u_i(t)|^2\le \epsilon_0^2,
\ \ {\rm{for}}\ \ i=1,2.
\end{equation}
Set $R_0=\min\{r_0, t_0^{\frac14}\}=t_0^\frac14>0$. Then  (\ref{serrin_small4}) implies
\begin{equation}\label{serrin_small5}
\max_{z=(x, t)\in\Omega\times [0, t_0]}\Big\|\nabla^2 u_i\Big\|_{L^\infty_tL^2_x(P_{R_0}(z)\cap (\Omega\times [0,t_0]))}\le \epsilon_0,
\ \ {\rm{for}}\ \ i=1,2.
\end{equation}
Hence $u_1$ and $u_2$ satisfy (\ref{e-cond-serrin}) of Theorem \ref{e-regularity-serrin} (with $p=2$ and $q=\infty$)
on $P_r(z)$, for any $z\in \Omega\times [0,t_0]$ and $r=\min\{R_0, d(x,\partial\Omega), t^\frac14\}>0$.
Hence by suitable scalings of the estimate of Theorem \ref{e-regularity-serrin}, we have
\begin{equation}\label{gradient_estimate70}
\max_{i,2} \Big|\nabla^m u_i(x,t)\Big|
\lesssim \epsilon_0 \left(\frac{1}{R_0^m}+\frac{1}{d^m(x,\partial\Omega)}+\frac{1}{t^{\frac{m}4}}\right)
\lesssim \epsilon_0\left(\frac{1}{d^m(x,\partial\Omega)}+\frac{1}{t^{\frac{m}4}}\right),  \ \forall\ m\ge 1,
\end{equation}
for any $(x,t)\in\Omega\times [0,t_0]$. Here we have used $R_0\ge t^\frac14$ in the last inequality.
Applying (\ref{gradient_estimate70}) and the proof of Theorem \ref{biharmonic-flow-unique1}, we can conclude
that $u_1\equiv u_2$ in $\Omega\times [0,t_0]$. Continuing this argument on the interval $[t_0,T]$ shows
$u_1\equiv u_2$ in $\Omega\times [0,T]$.  \qed

\medskip
\noindent{\bf Proof of Corollary \ref{4d_convex}}.  Let $\epsilon_2>0$ be given by Theorem \ref{e-regularity-serrin}.
Then (\ref{mono_decrease}) yields
\begin{equation}\label{small_serrin8}
\Big\|\nabla^2 u\Big\|_{L^\infty_tL^2_x(\Omega\times [0,\infty))}\le\epsilon_2.
\end{equation}
Hence by suitable scalings of the estimate of Theorem \ref{e-regularity-serrin}, we have $u\in C^\infty(\Omega\times (0, \infty), N)$ and there exists
$T_1>0$ such that
\begin{equation}\label{gradient_estimate90}
\Big|\nabla^m u(x,t)\Big|\lesssim \epsilon_2 \left(\frac{1}{d^m(x,\partial\Omega)}+\frac{1}{t^{\frac{m}4}}\right),
\ \forall\  m\ge 1,
\end{equation}
holds for all $x\in\Omega$ and $t\ge T_1$. Now we can apply the same arguments as in the proof of
Theorem \ref{biharmonic-flow-convex} and Corollary \ref{biharmonic-flow-unique2} to prove the
conclusions of Corollary \ref{4d_convex}. \qed

\section {Appendix: Higher order regularity}
\setcounter{equation}{0}
\setcounter{theorem}{0}

It is known, at least to experts, that higher order regularity holds for any H\"older continuous solution
to (\ref{biharmonic-flow}) of the heat flow of biharmonic maps . However,  we can't find a proof
in the literature. For the completeness, we will sketch  a proof here.

\begin{proposition}  For $0<\alpha<1$,  if $u\in W^{1,2}_2\cap C^\alpha (P_2, N)$ is a weak
solution of (\ref{biharmonic-flow}), then $u\in C^\infty(P_1,N)$, and
\begin{equation}
\Big\|\nabla^m u\Big\|_{C^0(P_1)}\lesssim \Big [u\Big]_{C^\alpha(P_2)}
+\Big\|u\Big\|_{L^2_tW^{2,2}_x(P_2)}, \ \forall\  m\ge 1. \label{grad_estimate50}
\end{equation}
\end{proposition}

\pf By  {\it Claim} 2 and {\it Claim} 3 in the proof of Theorem \ref{e-regularity-serrin},  it suffices to establish that
$\nabla^2 u\in M^{2,4-4\tilde\alpha}(P_\frac32)$ for some $\frac23<\tilde\alpha<1$,  and
\begin{equation}
\Big\|\nabla^2 u\Big\|_{M^{2,4-4\tilde\alpha}(P_\frac32)}\lesssim
\Big [u\Big]_{C^\alpha(P_2)}+\Big\|\nabla^2 u\Big\|_{L^2(P_2)}.
\label{morrey8}
\end{equation}
This will be achieved by the hole-filling type argument. For any fixed $z_0=(x_0,t_0)\in P_\frac32$ and
$0<r\le\frac14$, let $\phi\in C_0^\infty(\mathbb R^n)$ be a cut-off function of $B_r(x_0)$, i.e.,
$$0\le\phi\le 1, \ \phi\equiv 1 \ {\rm{in}}\ B_r(x_0), \ \phi\equiv 0 \ {\rm{outside}}\ B_{2r}(x_0),
\ |\nabla^m\phi|\le Cr^{-m}, \ \forall \  m\ge 1.$$
Set $\displaystyle c:=-\!\!\!\!\!\!\int_{P_r(z_0)} u\in \mathbb R^{L+1}$. Multiplying (\ref{biharmonic-flow}) by
$(u-c)\phi^4$ and integrating over $\mathbb R^n$, we obtain
\begin{eqnarray}\label{morrey9}
&&\frac{d}{dt}\int_{\mathbb R^n}|u-c|^2\phi^4
+2\int_{\mathbb R^n}\Delta (u-c)\cdot\Delta ((u-c)\phi^4)
=2\int_{\mathbb R^n}\mathcal N_{\hbox{bh}}[u]\cdot (u-c)\phi^4\nonumber\\
&\lesssim& \int_{\mathbb R^n} |\nabla^2 u|^2|u-c|\phi^4
+\int_{\mathbb R^n}|\nabla u||\nabla^2 u||\nabla((u-c)\phi^4)|.
\end{eqnarray}
For the second term in the left hand side of (\ref{morrey9}), we have
\begin{eqnarray}\label{morrey10}
&&2\int_{\mathbb R^n}\Delta (u-c)\cdot\Delta ((u-c)\phi^4)
=2\int_{\mathbb R^n}\nabla^2 (u-c)\cdot\nabla^2 ((u-c)\phi^4)\nonumber\\
&\ge& 2\int_{B_r(z_0)}|\nabla^2 u|^2-C\int_{\mathbb R^n}|u-c|^2(|\nabla^2\phi|^2+|\nabla\phi|^4)+
\phi^2|\nabla\phi|^2|\nabla u|^2.
\end{eqnarray}
Substituting (\ref{morrey10}) into (\ref{morrey9}) and integrating over $t\in [t_0-r^4, t_0]$, we obtain
\begin{eqnarray}\label{morrey11}
&&\int_{P_r(z_0)}|\nabla^2 u|^2
\leq \int_{B_{2r}(x_0)\times \{t_0-r^4\}}|u-c|^2
+\left(2^{-(n+4)}+C \hbox{osc}_{P_{2r}(z_0)}u\right) \int_{P_{2r}(z_0)}|\nabla^2 u|^2\nonumber\\
&&\qquad\qquad\quad\ \ \ \ \ \ +C r^{n}\left(\hbox{osc}_{P_{2r}(z_0)}u\right)^2+C\left[1+(\hbox{osc}_{P_{2r}(z_0)}u)^2\right] r^{-2}\int_{P_{2r}(z_0)}\phi^2|\nabla u|^2\nonumber\\
&&\qquad\qquad\quad\qquad+C\int_{P_{2r}(z_0)}|\nabla u|^4\phi^4\
\end{eqnarray}
By integration by parts and the H\"older inequality, we have
$$ r^{-2}\int_{P_{2r}(z_0)}\phi^2|\nabla u|^2\le
Cr^{-2}\left(\hbox{osc}_{P_{2r}(z_0)}u\right)\int_{P_{2r}(z_0)}|\nabla^2 u|+
Cr^n \left(\hbox{osc}_{P_{2r}(z_0)}u\right)^2,$$
and
$$C\int_{P_{2r}(z_0)}\phi^4|\nabla u|^4\le 2^{-(n+4)}\int_{P_{2r}(z_0)}|\nabla^2 u|^2+
Cr^n \left(\hbox{osc}_{P_{2r}(z_0)}u\right)^4+C\left(\hbox{osc}_{P_{2r}(z_0)}u\right)^2\int_{P_{2r}(z_0)}|\nabla^2 u|^2.$$
Putting these two inequalities into (\ref{morrey11}) and using
$\displaystyle\hbox{osc}_{P_{2r}(z_0)}u\le Cr^\alpha,$
we get
\begin{eqnarray}
&&\int_{P_r(z_0)}|\nabla^2 u|^2 \le\left(2^{-(n+3)}+Cr^\alpha\right)\int_{P_{2r}(z_0)}|\nabla^2 u|^2+Cr^{n+2\alpha}
+C(1+r^{2\alpha})r^{\alpha-2}\int_{P_{2r}(z_0)}|\nabla^2 u|\nonumber\\
&&\qquad\qquad\qquad\leq \left(2^{-(n+2)}+Cr^\alpha\right)\int_{P_{2r}(z_0)}|\nabla^2 u|^2+Cr^{n+2\alpha},
\label{morrey12}
\end{eqnarray}
where we have used the following inequality in the last step:
$$C(1+r^{2\alpha})r^{\alpha-2}\int_{P_{2r}(z_0)}|\nabla^2 u|
\leq 2^{-(n+3)}\int_{P_{2r}(z_0)}|\nabla^2 u|^2+C r^{n+2\alpha}.$$
Choosing $r>0$ so small that $Cr^\alpha\le 2^{-(n+3)}$,  we see that (\ref{morrey12}) implies
\begin{equation}\label{morrey13}
r^{-n}\int_{P_r(z_0)}|\nabla^2 u|^2\le \frac12 (2r)^{-n}\int_{P_{2r}(z_0)}|\nabla^2 u|^2 +Cr^{2\alpha}.
\end{equation}
It is clear that iterating (\ref{morrey13}) implies that there is $\alpha_0\in (0,1)$ such that
$\displaystyle\nabla^2 u\in M^{2,4-2\alpha_0}(P_\frac32)$ and
\begin{equation}\label{morrey14}
\Big\|\nabla^2 u\Big\|_{M^{2,4-2\alpha_0}(P_\frac32)}
\lesssim \Big [u\Big]_{C^\alpha(P_2)}+\Big\|\nabla^2 u\Big\|_{L^2(P_2)}.
\end{equation}
We can apply the estimate (\ref{morrey14}) and repeat the above argument to show that
$\displaystyle\nabla^2 u\in M^{2,4-4\alpha_0}(P_\frac32)$ and the estimate (\ref{morrey14}) holds with
$\alpha_0$ replaced by $2\alpha_0$. Repeating these argument again and again until
there exists $\tilde\alpha\in (\frac23,1)$ such that $\displaystyle\nabla^2 u\in M^{2,4-4\tilde\alpha}(P_\frac32)$
and the estimate (\ref{morrey8}) holds.  The remaining parts of the proof can be done by
following the same arguments as in
{\it Claim} 2 and {\it Claim} 3
of the proof of Theorem \ref{e-regularity-serrin}. This completes the proof.
\endpf

\bigskip
\noindent{\bf Acknowledgements}.
The authors are partially supported by NSF grant 1001115. The third author is also partially supported by 
a Simons Fellowship and NSFC grant 11128102.
The authors would like to thank the referee for many constructive suggestions that help to improve the presentation.

\end{document}